# INVASION PERCOLATION ON REGULAR TREES[1]

BY OMER ANGEL, JESSE GOODMAN, FRANK DEN HOLLANDER
AND GORDON SLADE

*University of British Columbia, University of Toronto and Leiden University*

We consider invasion percolation on a rooted regular tree. For the infinite cluster invaded from the root, we identify the scaling behavior of its $r$-point function for any $r \geq 2$ and of its volume both at a given height and below a given height. We find that while the power laws of the scaling are the same as for the incipient infinite cluster for ordinary percolation, the scaling functions differ. Thus, somewhat surprisingly, the two clusters behave differently; in fact, we prove that their laws are mutually singular. In addition, we derive scaling estimates for simple random walk on the cluster starting from the root. We show that the invasion percolation cluster is stochastically dominated by the incipient infinite cluster. Far above the root, the two clusters have the same law locally, but not globally.

A key ingredient in the proofs is an analysis of the forward maximal weights along the backbone of the invasion percolation cluster. These weights decay toward the critical value for ordinary percolation, but only slowly, and this slow decay causes the scaling behavior to differ from that of the incipient infinite cluster.

**1. Introduction and main results.**

1.1. *Motivation and background.* Invasion percolation is a stochastic growth model introduced by Wilkinson and Willemsen [17]. In its general setting, the edges of an infinite connected graph $\mathcal{G}$ are assigned i.i.d. uniform random variables on $(0, 1)$, called weights, a distinguished vertex $o$ is chosen, called the origin, and an infinite subgraph of $\mathcal{G}$ is grown inductively as follows. Define $I_0$ to be the vertex $o$. For $N \in \mathbb{N}_0$, given $I_N$, let $I_{N+1}$ be obtained by adjoining to $I_N$ the edge in its boundary with smallest

[1]Supported in part by NSERC of Canada.

weight. The *invasion percolation cluster* (IPC) is the random infinite subgraph $\bigcup_{N \in \mathbb{N}_0} I_N \subset \mathcal{G}$, which we denote by $C$. We will occasionally blur the distinction between $C$ as a graph and as a set of vertices.

Invasion percolation is closely related to critical percolation. Indeed, suppose $\mathcal{G}$ has a bond percolation threshold $p_c$ that lies strictly between 0 and 1, and color red those bonds (= edges) whose weight is at most $p_c$. Once a red bond is invaded, all other red bonds in its cluster will be invaded before the invasion process leaves the cluster. For $\mathcal{G} = \mathbb{Z}^d$, where critical clusters appear on all scales, we expect larger and larger critical clusters to be invaded, so that the invasion process spends a large proportion of its time in large critical clusters. A reflection of this is the fact, proved for $\mathcal{G} = \mathbb{Z}^d$ by Chayes, Chayes and Newman [5] and extended to much more general graphs by Häggström, Peres and Schonmann [6], that the number of bonds in $C$ with weight above $p_c + \varepsilon$ is almost surely finite, for all $\varepsilon > 0$. When $\mathcal{G}$ is a regular tree, this fact is easy to prove: For any $p > p_c$, whenever an edge is invaded with weight above $p$, there is an independent positive probability of encountering an infinite cluster consisting of edges of weight at most $p$, and never again invading an edge of weight above $p$. Therefore, the number of invaded edges above $p$ is finite. The fact that invasion percolation is driven by the critical parameter $p_c$, even though there is no parameter specification in its definition, makes it a prime example of *self-organized criticality*.

Another reflection of the relation to critical percolation has been obtained by Járai [11], who showed for $\mathbb{Z}^2$ that the probability of an event $E$ under the *incipient infinite cluster* (IIC) measure (constructed by Kesten [12]) is identical to the probability of the translation of $E$ to $x \in \mathbb{Z}^2$ under the IPC measure, conditional on $x$ being invaded and in the limit as $\|x\| \to \infty$. It is tempting to take this a step further and *conjecture* that the scaling limit of invasion percolation on $\mathbb{Z}^d$ when $d > 6$ is the canonical measure of super-Brownian motion conditioned to survive forever (see van der Hofstad [8], Conjecture 6.1). Indeed, such a result was proved for the IIC of spread-out (= long-range) oriented percolation on $\mathbb{Z}^d \times \mathbb{N}_0$ when $d > 4$ in van der Hofstad, den Hollander and Slade [9], and van der Hofstad [8], and presumably it holds for the IIC of unoriented percolation on $\mathbb{Z}^d$ when $d > 6$ as well.

Invasion percolation on a regular tree was studied by Nickel and Wilkinson [16]. They computed the probability generating function for the height and weight of the bond added to $I_N$ to form $I_{N+1}$. They looked, in particular, at the expected number of vertices in $I_N$ at level $t\sqrt{N}$, for $t \in [0, \infty]$ fixed and $N \to \infty$, and found that this expectation is described by the same power law as in critical percolation, but has a different dependence on $t$ (i.e., has a different shape function). They refer to this discrepancy as the "paradox of invasion percolation." Their analysis does not apply directly to the infinite IPC, so it does not allow for a direct comparison with the IIC. It does suggest though that *the* IPC *has a different scaling limit than the* IIC.

Let $\mathcal{T}_\sigma$ denote the rooted regular tree with forward degree $\sigma \geq 2$ (i.e., all vertices have degree $\sigma + 1$, except the root $o$, which has degree $\sigma$). In the present paper, we study the IPC on $\mathcal{T}_\sigma$ (see Figure 1 for a simulation), and show that indeed it does *not* have the same scaling limit as the IIC. Furthermore, we show that the laws of the IPC and the IIC are mutually singular. There is no reason to believe that this discrepancy will disappear for other graphs, such as $\mathbb{Z}^d$, and so the conjecture raised in [8] must be expected to be *false*.

Central to our analysis is a representation of $C$ as *an infinite backbone* (an infinite self-avoiding path rising from the root) *from which emerge branches having the same distribution as subcritical percolation clusters.* The percolation parameter value of these subcritical branches depends on a process we call the *forward maximal weight process* along the backbone. We analyze this process in detail, and prove, in particular, that as $k \to \infty$ the maximum weight of a bond on the backbone above height $k$ is asymptotically $p_c(1 + Z/k)$, where $Z$ is an exponential random variable with mean 1. This quantifies the rate at which maximal bond weights approach $p_c$ as the invasion proceeds. It is through an understanding of this process that the "paradox of invasion percolation" can be resolved, both qualitatively and quantitatively.

It is interesting to compare the above slow decay with the inhomogeneous model of Chayes, Chayes and Durrett [4], in which the percolation parameter $p$ depends on $x \in \mathbb{Z}^d$ and scales like $p_c + \|x\|^{-(\epsilon+1/\nu)}$, where $\nu$ is the critical exponent for the correlation length. It is proved in [4] that for $\mathbb{Z}^2$ (and conjectured for $\mathbb{Z}^d$ for $d > 2$) that when $\varepsilon < 0$ the origin has a positive probability of being in an infinite cluster, but not when $\varepsilon > 0$. For invasion percolation on a tree, the weight $p_c(1 + Z/k)$ corresponds to the boundary value $\varepsilon = 0$ (we use graph distance on the tree), but with a random coefficient $Z$. Invasion percolation, therefore, corresponds in some sense to the critical case of the inhomogeneous model.

From our analysis of the forward maximal weight process along the backbone of invasion percolation on a tree, we are able to compute the scaling of all the $r$-point functions of $C$, and of the size of $C$ both at a given height and below a given height. The scaling limits are independent of $\sigma$ apart from a simple overall factor. Each of these quantities scales according to the same powers laws as their counterparts for the IIC, but with different scaling functions. The Hausdorff dimension of both clusters is 4. Moreover, we apply results of Barlow, Járai, Kumagai and Slade [1] to prove scaling estimates for simple random walk on $C$ starting from $o$. These estimates establish that $C$ has spectral dimension $\frac{4}{3}$, which is the same as for the IIC (see also Kesten [13], and Barlow and Kumagai [2]).

It would be of interest to extend our results to invasion percolation on $\mathbb{Z}^d$ when $d > 6$ in the unoriented setting and on $\mathbb{Z}^d \times \mathbb{N}_0$ when $d > 4$ in the oriented setting, where lace expansion methods could be tried. However, it seems a challenging problem to carry over the expansion methods developed in Hara and Slade [7], van der Hofstad and Slade [10], and Nguyen and Yang [15], since invasion percolation lacks bond independence and uses supercritical bonds. An additional motivation for the problem on $\mathbb{Z}^d$ is the following observation of Newman and Stein [14]: if the probability that $x \in C$ scales like $\|x\|^{4-d}$, then this has consequences for the number of ground states of a spin glass model when $d > 8$.

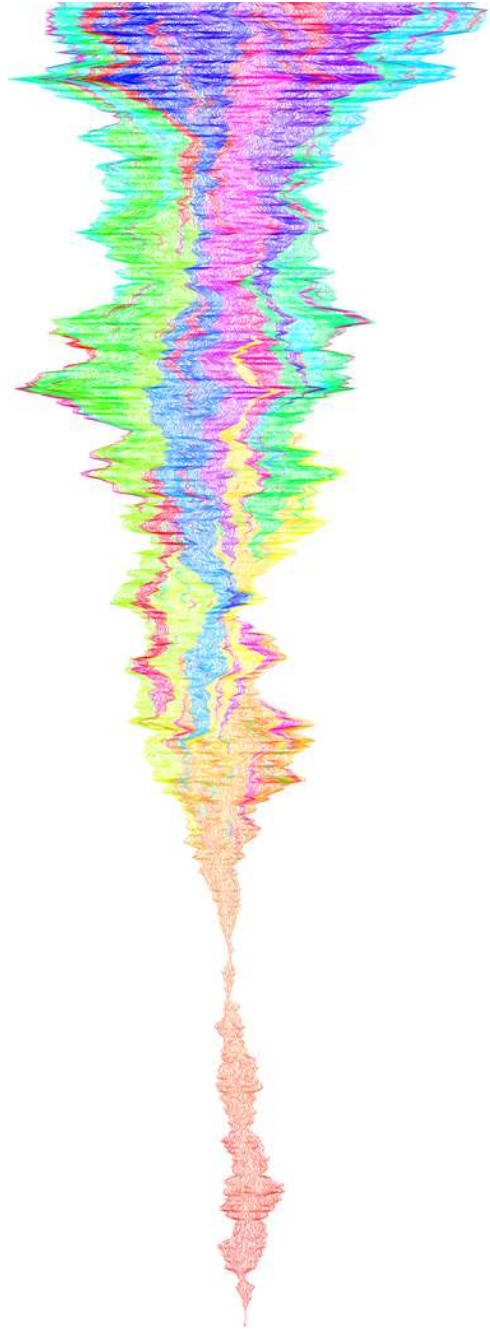

FIG. 1. *Simulation of invasion percolation on the binary tree up to height 500. The hue of the $i$th added edge is $i/M$, with $M$ the number of edges in the figure. The color sequence is red, orange, yellow, green, cyan, blue, purple and red. The last edge is almost as red as the first.*

We begin in Section 1.2 with a review of the IIC on $\mathcal{T}_\sigma$, for later comparison with our results for the IPC, which are stated in Section 1.3. Section 1.4 outlines the rest of the paper.

Before discussing the IIC, we introduce some notation. We denote the height of a vertex $v \in \mathcal{T}_\sigma$ by $\|v\|$; this is its graph distance from $o$ in $\mathcal{T}_\sigma$. We write $\mathbb{P}_p$ for the law of independent bond percolation with parameter $p$, $\mathbb{P}_\infty$ for the law of the IIC of independent bond percolation, and $\mathbb{P}$ for the law of the IPC.

1.2. *The incipient infinite cluster.* The IIC on a tree is discussed in detail in Kesten [13] and in Barlow and Kumagai [2]. It is constructed by conditioning a critical branching process to survive until height $n$, and then letting $n \to \infty$. In our case, the branching process has a binomial offspring distribution with parameters $(\sigma, 1/\sigma)$. We summarize some elementary properties of the IIC in this section. To keep our exposition self-contained, we provide quick indications of proofs of these properties in Section 9.

On $\mathcal{T}_\sigma$, the IIC can be viewed as consisting of *an infinite backbone adorned with branches at each vertex that are independent critical percolation clusters in each direction away from the backbone.* We write $C_\infty$ to denote the IIC. This is an infinite random subgraph of $\mathcal{T}_\sigma$, but it will be convenient to think of $C_\infty$ as a set of vertices.

Fix $r \geq 2$. Pick $r-1$ vertices $\vec{x} = (x_1, \ldots, x_{r-1})$ in $\mathcal{T}_\sigma \setminus \{o\}$ such that no $x_i$ lies on the path from $o$ to any $x_j$ ($j \neq i$). Let $\mathcal{S}(\vec{x})$ denote the subtree of $\mathcal{T}_\sigma$ obtained by connecting the vertices in $\vec{x}$ to $o$. Call this the *spanning tree* of $o$ and $\vec{x}$. Let $N$ denote the number of edges in $\mathcal{S}(\vec{x})$. Write $\vec{x} \in C_\infty$ for the event that all vertices in $\vec{x}$ lie in $C_\infty$, which is the same as the event that $\mathcal{S}(\vec{x}) \subset C_\infty$. The $r$-point function is the probability $\mathbb{P}_\infty(\vec{x} \in C_\infty)$ (with $o$ the $r$th point). Let $\partial \mathcal{S}(\vec{x})$ denote the external boundary of $\mathcal{S}(\vec{x})$; this is the set of vertices in $\mathcal{T}_\sigma \setminus \mathcal{S}(\vec{x})$ whose parent is a vertex in $\mathcal{S}(\vec{x})$. The cardinality of $\partial \mathcal{S}(\vec{x})$ is $N(\sigma - 1) + \sigma$. For $y \in \partial \mathcal{S}(\vec{x})$, let $B_y$ denote the event that $y$ is in the backbone, that is, $y$ is the first vertex in the backbone after it emerges from $\mathcal{S}(\vec{x})$. Then

$$\begin{aligned}\sigma^{N+1}\mathbb{P}_\infty(\vec{x} \in C_\infty) &= N(\sigma - 1) + \sigma, \\ \mathbb{P}_\infty(B_y \mid \vec{x} \in C_\infty) &= \frac{1}{N(\sigma - 1) + \sigma}, \qquad y \in \partial \mathcal{S}(\vec{x}).\end{aligned} \tag{1.1}$$

The first line of (1.1) gives a simple formula for the $r$-point function of the IIC, in which only the size of $\mathcal{S}(\vec{x})$ is relevant, not its geometry. The second line shows that the backbone emerges uniformly from $\mathcal{S}(\vec{x})$.

Let

$$
\begin{aligned}
C_\infty[n] &= \sharp\{x \in C_\infty : \|x\| = n\}, \\
C_\infty[0,n] &= \sharp\{x \in C_\infty : 0 \leq \|x\| \leq n\}, \qquad n \in \mathbb{N}_0,
\end{aligned}
\tag{1.2}
$$

and abbreviate

$$
\rho = \rho(\sigma) = \frac{\sigma - 1}{2\sigma}.
\tag{1.3}
$$

Then, under the law $\mathbb{P}_\infty$,

$$
\frac{1}{\rho n} C_\infty[n] \Longrightarrow \Gamma_\infty, \qquad \frac{1}{\rho n^2} C_\infty[0,n] \Longrightarrow \widehat{\Gamma}_\infty, \qquad n \to \infty,
\tag{1.4}
$$

where $\Longrightarrow$ denotes convergence in distribution, and $\Gamma_\infty, \widehat{\Gamma}_\infty$ are random variables with Laplace transforms

$$
\mathbb{E}_\infty(e^{-\tau \Gamma_\infty}) = (1+\tau)^{-2}, \qquad \mathbb{E}_\infty(e^{-\tau \widehat{\Gamma}_\infty}) = [\cosh(\sqrt{\tau})]^{-2}, \qquad \tau \geq 0.
\tag{1.5}
$$

$\Gamma_\infty$ is the size biased exponential with parameter 1, that is, the distribution with density $xe^{-x}$, $x \geq 0$. It is straightforward to compute the moments:

$$
\mathbb{E}_\infty(\Gamma_\infty) = 2, \qquad \mathbb{E}_\infty(\Gamma_\infty^2) = 6, \qquad \mathbb{E}_\infty(\widehat{\Gamma}_\infty) = 1, \qquad \mathbb{E}_\infty(\widehat{\Gamma}_\infty^2) = \tfrac{4}{3}.
\tag{1.6}
$$

1.3. *Main results.* This section contains our main results for the scaling behavior of $C$ under the law $\mathbb{P}$, listed in Sections 1.3.1–1.3.5.

It is easy to see that, under the law $\mathbb{P}$, $C$ has almost surely a *single* backbone. Indeed, suppose that with positive $\mathbb{P}$-probability there is a vertex in $C$ from which there are two disjoint paths to infinity. Conditioned on this event, let $M_1$ and $M_2$ denote the maximal weights along these paths. It is not possible that $M_1 > M_2$, because the entire infinite second branch would be invaded before the edge carrying the weight $M_1$; $M_2 > M_1$ is ruled out for the same reason. However, $M_1 = M_2$ has probability zero, because the distribution of the weights is continuous.

1.3.1. *Stochastic domination and local behavior.* The following two theorems will be proved in Section 2. The first theorem is part of a deeper structural representation of the IPC, which is described in Section 2.1 and which is the key to all our scaling results.

THEOREM 1.1. *The* IIC *stochastically dominates the* IPC, *that is, there exists a coupling of $C_\infty$ and $C$ such that $C_\infty \supset C$ with probability 1.*

THEOREM 1.2. *Let $\mathcal{T}_\sigma^*$ denote the rooted regular tree in which all vertices (including the root) have degree $\sigma + 1$. Let $E$ be a cylinder event on $\mathcal{T}_\sigma^*$ (i.e., an event that depend on the status of only finitely many bonds), and suppose that $E$ is invariant under the automorphisms of $\mathcal{T}_\sigma^*$. Then*

$$\lim_{\|x\| \to \infty} \mathbb{P}(\tau_x E \mid x \in C) = \mathbb{P}_\infty^*(E), \tag{1.7}$$

*where $\tau_x$ denotes the shift by $x$, and $\mathbb{P}_\infty^*$ denotes the* IIC *on $\mathcal{T}_\sigma^*$.*

The symmetry assumption on $E$ in Theorem 1.2 is necessary because the unique path in the tree from $o$ to $x$ must be invaded when $x \in C$, whereas $\mathbb{P}_\infty^*$ has no such preferred path. Theorem 1.2 shows that $C$ and $C_\infty$ are the same *locally* far above $o$. Comparing the results in Sections 1.3.2 – 1.3.3 below with the analogous results for the IIC show that *globally* they are different.

Járai [11] proves additional statements in the spirit of Theorem 1.2 for invasion percolation on $\mathbb{Z}^2$. We expect that similar statements can be proved also for the tree, but we do not pursue these here.

1.3.2. *The $r$-point function.* For $r \geq 2$, the invasion percolation $r$-point function is the probability $\mathbb{P}(x_1, \ldots, x_{r-1} \in C)$, which we write simply as $\mathbb{P}(\vec{x} \in C)$ with $\vec{x} = (x_1, \ldots, x_{r-1})$. We can and do assume that no $x_i$ lies on the path from $o$ to any $x_j$ ($j \neq i$), since any such $x_i$ is automatically invaded when $x_j$ is.

To state our result for the asymptotics of the $r$-point function, some more terminology is required. We recall the definition of $\mathcal{S}(\vec{x})$, $\partial \mathcal{S}(\vec{x})$, $N$ and $B_y$ given in Section 1.2. Let $\mathcal{N}(\vec{x})$ denote the set of *nodes* of $\mathcal{S}(\vec{x})$; this is the set consisting of $o$, the $r-1$ vertices in $\vec{x}$ and any additional vertices where $\mathcal{S}(\vec{x})$ branches. For $v \in \mathcal{N}(\vec{x}) \setminus \{o\}$, write $v_-$ to denote the node immediately below $v$, and $n_v$ to denote the number of edges in the segment of $\mathcal{S}(\vec{x})$ between $v_-$ and $v$. We write $w < v$ when $w$ is a node below $v$. For $w, v \in \mathcal{N}(\vec{x})$ with $w < v$, let $M_w^v$ denote the number of edges in the subtree obtained from $\mathcal{S}(\vec{x})$ by *deleting everything above $w$ in the direction of $v$*. (See Figure 2 for an illustration.)

Given $y \in \partial \mathcal{S}(\vec{x})$, let $v$ be the first node above or equal to the parent of $y$, and let $k$ be the distance from $v_-$ to the parent of $y$. Note that $v$ and $k$ depend on $y$, but we will *not* make this explicit in our notation.

Theorem 1.3 and Corollary 1.4, which will be proved in Section 4, describe a scaling limit in which the *lengths* of all the segments of $\mathcal{S}(\vec{x})$ tend to infinity while the *geometry* of $\mathcal{S}(\vec{x})$ stays the same. More precisely, given $t_v \in (0,1)$ for each $v \in \mathcal{N}(\vec{x}) \setminus \{o\}$, with $\sum_{v \in \mathcal{N}(\vec{x}) \setminus \{o\}} t_v = 1$, we assume that

$$\frac{n_v}{N} \to t_v, v \in \mathcal{N}(\vec{x}) \setminus \{o\} \qquad \text{as } N \to \infty \tag{1.8}$$

and, given $s \in [0, t_v]$, that

$$\frac{k}{N} \to s \qquad \text{as } N \to \infty, \tag{1.9}$$

with $k$ and $v$ related to $y$ as described above. We write $\sharp\lim_{N\to\infty}$ to denote the limit in (1.8)–(1.9). Furthermore, we define

$$\sharp\lim_{N\to\infty} \frac{M_w^v}{N} = m_w^v, \qquad w, v \in \mathcal{N}(\vec{x})\setminus\{o\}, w < v. \tag{1.10}$$

In the scaling limit, we may associate with $\mathcal{S}(\vec{x})$ and $\mathcal{N}(\vec{x})$ a *scaled spanning tree* $\mathcal{S}$ with nodes $\mathcal{N}$. The segments of this tree are labeled by $\mathcal{N}\setminus\{o\}$ and are continuous line pieces with lengths $t_v$, $v \in \mathcal{N}\setminus\{o\}$. The backbone emerges at height $s$ above the bottom of segment $v$.

THEOREM 1.3. *Let $r \geq 2$. Suppose that $\mathcal{S}$ does not branch at $o$ (i.e., $o$ has degree 1 in $\mathcal{S}$). Then*

$$\sharp\lim_{N\to\infty} \sigma^{N+1}\mathbb{P}(\vec{x} \in C, B_y) = (s + m_{v_-}^v)\pi_v, \qquad y \in \partial\mathcal{S}(\vec{x}), \tag{1.11}$$

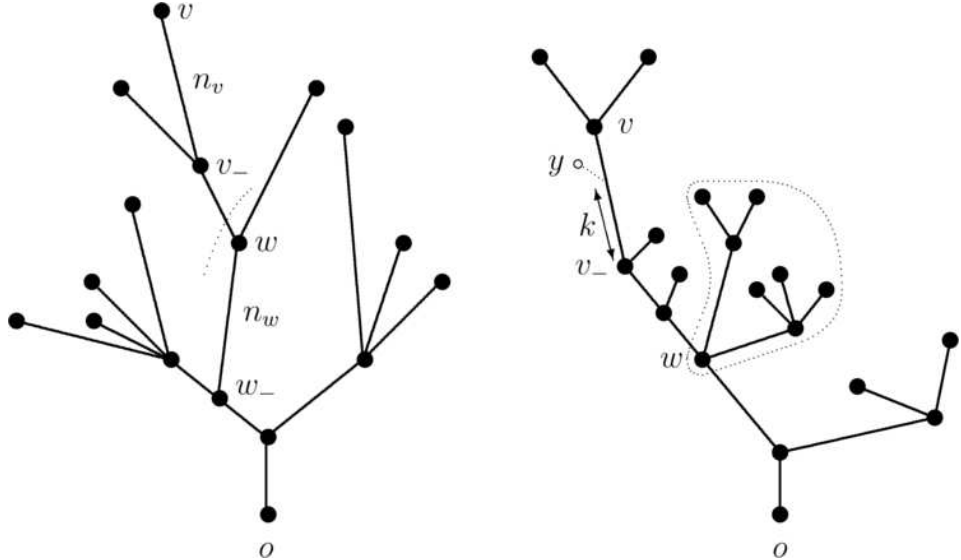

FIG. 2. *The illustration at left shows a spanning tree $\mathcal{S}(\vec{x})$ for $r = 11$. The dots are the nodes in $\mathcal{N}(\vec{x})$. The dots at the leaves are the vertices in $\vec{x}$. The dotted line indicates the cut that deletes everything above $w$ in the direction of $v$; $M_w^v$ is the number of edges left after the cut. The illustration at right, for $r = 12$, shows the relation between $y, v, v_-, k$, and the dotted line isolates the edges contributing to $N_w^v$ defined in Section 4.*

*where*

$$\pi_v = \prod_{\substack{w \in \mathcal{N} \\ o < w < v}} \frac{t_w + m_{w_-}^v}{m_w^v} \tag{1.12}$$

*with the convention that the empty product is* 1.

Note that in the right-hand side of (1.11) the dependence on $s$ is linear, and that $\pi_v$ and $m_{v_-}^v$ are simple functionals of the geometry of the scaled spanning tree $\mathcal{S}$. Further note that $\pi_v$ is a product of ratios that take values in $(0, 1)$.

By summing (1.11) over $y \in \partial \mathcal{S}(\vec{x})$, which amounts to summing first over $0 < k \leq n_v$ and then over $v \in \mathcal{N}(\vec{x}) \setminus \{o\}$, we will derive the asymptotics for the $r$-point function.

COROLLARY 1.4. *Let $r \geq 2$. Suppose that $\mathcal{S}$ does not branch at $o$. Then*

$$\sharp \lim_{N \to \infty} \frac{1}{(\sigma - 1)N} \sigma^{N+1} \mathbb{P}(\vec{x} \in C) = \sum_{v \in \mathcal{N} \setminus \{o\}} \left( \frac{1}{2} t_v^2 + t_v m_{v_-}^v \right) \pi_v. \tag{1.13}$$

By combining (1.11)–(1.13), we obtain the distribution for the vertex where the backbone emerges from $\mathcal{S}(\vec{x})$, conditional on $\mathcal{S}(\vec{x})$ being invaded:

$$\begin{aligned}
\sharp \lim_{N \to \infty} (\sigma - 1) N \mathbb{P}(B_y \mid \vec{x} \in C) \\
= \frac{(s + m_{v_-}^v) \pi_v}{\sum_{u \in \mathcal{N} \setminus \{o\}} ((1/2) t_u^2 + t_u m_{u_-}^u) \pi_u}, \qquad y \in \partial \mathcal{S}(\vec{x}).
\end{aligned} \tag{1.14}$$

The restriction in Theorem 1.3 and Corollary 1.4 that $\mathcal{S}$ does not branch at $o$ is essential. We will see in Section 4 that when $\mathcal{S}$ branches at $o$ the limit in (1.11) is zero for all $y \in \mathcal{S}(\vec{x})$, that is, diagrams branching at the bottom are of higher order.

The following two examples illustrate (1.13)–(1.14):

*Two-point function*: For $r = 2$, $\mathcal{S}(\vec{x})$ consists of $o$ and a single vertex $x_1$ at height $n_1 = N$. See Figure 3. In this case, $m_o^1 = 0$ and $\pi_1 = 1$, and therefore

$$\begin{aligned}
\sharp \lim_{N \to \infty} \frac{1}{(\sigma - 1)N} \sigma^{N+1} \mathbb{P}(x_1 \in C) = \frac{1}{2}, \\
\sharp \lim_{N \to \infty} (\sigma - 1) N \mathbb{P}(B_y \mid x_1 \in C) = 2s, \qquad y \in \partial \mathcal{S}(x_1).
\end{aligned} \tag{1.15}$$

The first formula in (1.15) also follows directly from the results of Nickel and Wilkinson [16]. The second formula in (1.15) shows that the backbone branches off the path from $o$ to $x_1$ with an asymptotically *linear* density. This should be contrasted with the *constant* density in (1.1) for the IIC. In particular, the backbone for invasion percolation is more likely to branch off later than earlier. The reason for this will be discussed at the end of Section 2.1.

*Three-point function*: For $r = 3$, $\mathcal{S}(\vec{x})$ consists of the nodes $o$, $x_*$ at height $n_*$, and $x_1, x_2$ at heights $n_1, n_2$ above $x_*$. See Figure 3. By definition, $m_*^1 = t_* + t_2$, $m_*^2 = t_* + t_1$, $\pi_* = 1$, $\pi_1 = t_*/(t_* + t_2)$, and $\pi_2 = t_*/(t_* + t_1)$. Let

$$(1.16) \qquad u(t_*, t_1, t_2) = \frac{1}{2}\left(1 + \frac{t_1}{t_* + t_2} + \frac{t_2}{t_* + t_1}\right).$$

Then, after some arithmetic, we find that

$$(1.17) \qquad \sharp \lim_{N \to \infty} \frac{1}{(\sigma - 1)N} \sigma^{N+1} \mathbb{P}(x_1, x_2 \in C) = t_* u(t_*, t_1, t_2)$$

and

$$(1.18) \quad \begin{aligned} &\sharp \lim_{N \to \infty} (\sigma - 1) N \mathbb{P}(B_y \mid x_1, x_2 \in C) \\ &= \frac{1}{u(t_*, t_1, t_2)} \times \begin{cases} \dfrac{1}{t_*} s_*, & y \in \partial \mathcal{S}_*(\vec{x}), \\ 1 + \dfrac{1}{t_* + t_2} s_1, & y \in \partial \mathcal{S}_1(\vec{x}), \\ 1 + \dfrac{1}{t_* + t_1} s_2, & y \in \partial \mathcal{S}_2(\vec{x}), \end{cases} \end{aligned}$$

where $\partial \mathcal{S}_*(\vec{x}), \partial \mathcal{S}_1(\vec{x}), \partial \mathcal{S}_2(\vec{x})$ denote the external boundaries of the respective segments of $\mathcal{S}(\vec{x})$. Note that the right-hand side of (1.18) is a density on the scaled spanning tree $\mathcal{S}$ that is *linearly increasing on each segment*, and is *continuous at the nodes*.

A similar picture follows from (1.14) for all $r \geq 2$. The linear slope depends on the structure of the subtree obtained by cutting off everything above the segment, and decreases when moving upward in the tree. This is in sharp contrast with the uniform distribution for the IIC in (1.1), and shows that the scaling limits of the IPC and the IIC are different.

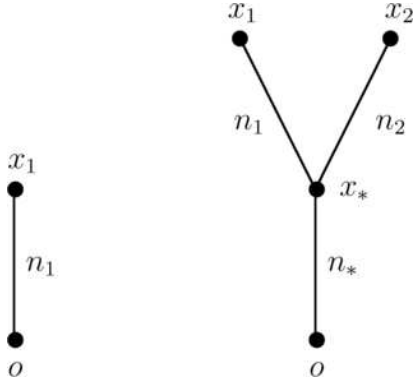

FIG. 3. *Spanning trees for $r = 2$ and $r = 3$.*

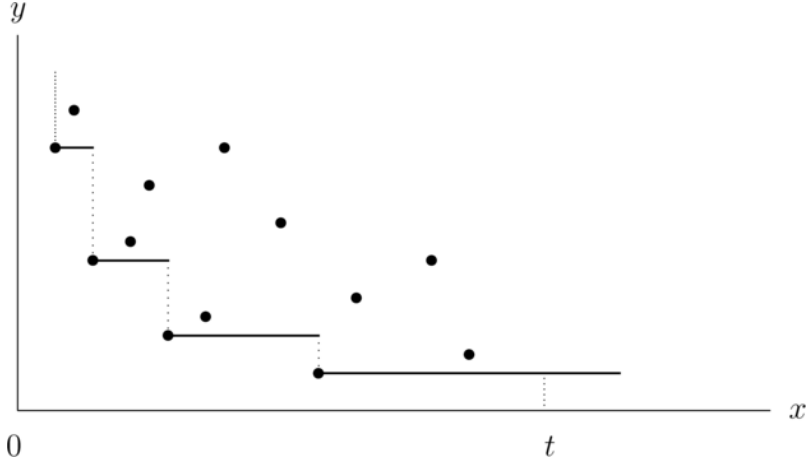

FIG. 4. *Sketch of the graph of $L(t)$ versus $t$. The dots are the points in $\mathcal{P}$.*

1.3.3. *Cluster size asymptotics.* Let $\mathcal{P}$ denote the Poisson point process on the positive quadrant with intensity 1. Write $\mathbb{P}_\mathcal{P}$ to denote its law. Let $L\colon (0,\infty) \to (0,\infty)$ denote its *lower envelope*, defined by

$$(1.19) \qquad L(t) = \min\{y > 0 : (x,y) \in \mathcal{P} \text{ for some } x \leq t\}, \qquad t > 0.$$

See Figure 4 for an illustration. This is a cadlag process, piecewise constant and nonincreasing, with $\lim_{t\downarrow 0} L(t) = \infty$ and $\lim_{t\to\infty} L(t) = 0$, $\mathbb{P}_\mathcal{P}$-a.s. In Section 3.2, we will compute its multivariate Laplace transform.

As in (1.2), let $C[n]$ denote the number of vertices in $C$ at height $n$, and let $C[0,n] = \sum_{m=0}^n C[m]$ denote the number of vertices up to height $n$. Recall from (1.3) that $\rho = (\sigma - 1)/2\sigma$.

THEOREM 1.5. *Let $\Gamma_n = \frac{1}{\rho n} C[n]$. Under the law $\mathbb{P}$, $\Gamma_n \Longrightarrow \Gamma$ as $n \to \infty$, where $\Gamma$ is the random variable with Laplace transform*

$$(1.20) \qquad \mathbb{E}(e^{-\tau\Gamma}) = \mathbb{E}_\mathcal{P}(e^{-S(\tau,L)}), \qquad \tau \geq 0,$$

*with*

$$(1.21) \qquad S(\tau, L) = 2\tau \int_0^1 dt \frac{L(t) e^{-(1-t)L(t)}}{L(t) + \tau[1 - e^{-(1-t)L(t)}]}.$$

We will show in Section 5 that

$$(1.22) \qquad \lim_{n\to\infty} \mathbb{E}(\Gamma_n) = \mathbb{E}(\Gamma) = 1, \qquad \lim_{n\to\infty} \mathbb{E}(\Gamma_n^2) = \mathbb{E}(\Gamma^2) = \tfrac{5}{3}.$$

THEOREM 1.6. *Let $\widehat{\Gamma}_n = \frac{1}{\rho n^2} C[0,n]$. Under the law $\mathbb{P}$, $\widehat{\Gamma}_n \Longrightarrow \widehat{\Gamma}$ as $n \to \infty$, where $\widehat{\Gamma}$ is the random variable with Laplace transform*

$$\mathbb{E}(e^{-\tau \widehat{\Gamma}}) = \mathbb{E}_{\mathcal{P}}(e^{-\widehat{S}(\tau, L)}), \qquad \tau \geq 0, \tag{1.23}$$

*with*

$$\widehat{S}(\tau, L) = 4\tau \int_0^1 \frac{dt}{L(t) + \kappa(\tau, t) \coth[(1/2)(1-t)\kappa(\tau, t)]}, \tag{1.24}$$

*and $\kappa(\tau, t) = \sqrt{4\tau + L(t)^2}$.*

We will show in Section 6 that

$$\lim_{n \to \infty} \mathbb{E}(\widehat{\Gamma}_n) = \mathbb{E}(\widehat{\Gamma}) = \tfrac{1}{2}, \qquad \lim_{n \to \infty} \mathbb{E}(\widehat{\Gamma}_n^2) = \mathbb{E}(\widehat{\Gamma}^2) = \tfrac{25}{72}. \tag{1.25}$$

We see no way to evaluate the expectations in (1.20) and (1.23) in closed form, despite our knowledge of the multivariate Laplace transform of the $L$-process. Theorems 1.5–1.6, in addition to showing that the two scaling limits exist, exhibit the underlying complexity of the IPC and underline the key role that is played by the $L$-process. We will see in Section 9 that by setting $L \equiv 0$, we recover the expressions for the IIC in (1.5).

The laws of $\Gamma$ and $\widehat{\Gamma}$ are not the same as their IIC counterparts $\Gamma_\infty$ and $\widehat{\Gamma}_\infty$, as is immediate from a comparison of (1.22) and (1.25) with (1.6). The power law scalings of $C[n]$ and $C[0,n]$ in Theorems 1.5–1.6 are, however, the same linear and quadratic scalings as for the IIC. In particular, Theorem 1.6 is a statement that the Hausdorff dimension of the IPC is 4, as it is for the IIC. (For this, we imagine that paths in the IPC are embedded in $\mathbb{Z}^d$ as random walk paths, with the root mapped to the origin, so that the on the order of $n^2 = r^4$ vertices in the IPC below level $n = r^2$ will be within distance $r$ of the origin.) Comparing the values of the first and second moments of $\widehat{\Gamma}$ and $\widehat{\Gamma}_\infty$, we see that the IPC has half the size of the IIC on average, while the ratio of the variance of the size of the IPC to the square of its mean is $\tfrac{7}{18}$, compared to $\tfrac{1}{3}$ for the IIC. The relatively larger fluctuation for the IPC is due to the randomness of the weights on the backbone; this will be discussed further in Section 2.1.

The scaling of the first and second moments of $C[n]$ and $C[0,n]$ implied by (1.22) and (1.25) can also be deduced directly from the scaling of the 2-point and the 3-point function [recall (1.15) and (1.17)]. In the same manner we can deduce that

$$\lim_{\substack{n_1, n_2 \to \infty \\ n_1/n_2 \to a}} \mathbb{E}(\Gamma_{n_1} \Gamma_{n_2}) = 1 + \tfrac{1}{3} a(1+a), \qquad a \in [0,1], \tag{1.26}$$

as we will show in Section 5.3. It would be interesting to study $(\Gamma_n)_{n \in \mathbb{N}}$ as a process, but we do not pursue this here.

1.3.4. *Mutual singularity of IPC and IIC.* The following theorem is essentially a consequence of Theorem 1.5. It shows a dramatic manifestation of the difference between the IPC and the IIC.

THEOREM 1.7. *The laws of* IPC *and* IIC *are mutually singular.*

1.3.5. *Simple random walk on the invasion percolation cluster.* Given $C$, let $\mu_y$ denote the degree in $C$ (both forward and backward) of a vertex $y \in C$. Consider the discrete-time simple random walk $X = (X_k)_{k \in \mathbb{N}_0}$ on $C$ that starts at $X_0 = x$ and makes transitions from $y$ in $C$ to any neighbor of $y$ in $C$ with probability $1/\mu_y$. Denote the law of this random walk given $C$ by $P_C^x$, with corresponding expectation $E_C^x$. We will consider three quantities:

$$(1.27) \qquad R_k = \{X_0, \ldots, X_k\},$$

the range of $X$ up to time $k$, with cardinality $|R_k|$; the $k$-step transition kernel

$$(1.28) \qquad p_k^C(x, y) = \frac{1}{\mu_y} P_C(X_k = y \mid X_0 = x),$$

which satisfies the reversibility relation $p_k^C(x, y) = p_k^C(y, x)$; the first exit time above height $n$, $T_n = \min\{k \geq 0 : \|X_k\| = n\}$. The following theorem provides power laws for these three quantities.

THEOREM 1.8. *There is a set $\Omega_0$ of configurations of the* IPC *with $\mathbb{P}(\Omega_0) = 1$, and positive constants $\alpha_1, \alpha_2$, such that for each configuration $C \in \Omega_0$ and for each $x \in C$, the simple random walk on $C$ obeys the following:*

(a)

$$(1.29) \qquad \lim_{k \to \infty} \frac{\log |R_k|}{\log k} = \frac{2}{3}, \qquad P_C^x\text{-}a.s.$$

(b) *There exists $K_x(C) < \infty$ such that*

$$(1.30) \quad (\log k)^{-\alpha_1} k^{-2/3} \leq p_{2k}^C(x, x) \leq (\log k)^{\alpha_1} k^{-2/3} \qquad \forall k \geq K_x(C).$$

(c) *There exists $N_x(C) < \infty$ such that*

$$(1.31) \qquad (\log n)^{-\alpha_2} n^3 \leq E_C^x(T_n) \leq (\log n)^{\alpha_2} n^3 \qquad \forall n \geq N_x(C).$$

The results in Theorem 1.8 are similar to the behavior of simple random walk on the IIC; see Barlow, Járai, Kumagai and Slade [1], Barlow and Kumagai [2], Kesten [13]. The *spectral dimension* $d_s$ of $C$ can be defined by

$$(1.32) \qquad d_s = -2 \lim_{k \to \infty} \frac{\log p_{2k}^C(o, o)}{\log k}.$$

From (1.30) we see that $d_s = \frac{4}{3}$. For additional statements concerning the height $\|X_n\|$ after $n$ steps, see [1].

With the help of results from [2], it is shown in [1], Example 1.9(ii), that (1.29)–(1.31) hold for simple random walk on *any* random subtree of the IIC for $\mathcal{T}_\sigma$ such that the expectation of $1/C[0,n]$ is bounded above by a multiple of $1/n^2$. In view of Theorem 1.1, to prove Theorem 1.8, it therefore, suffices to prove the following uniform bound, which will be done in Section 8.

THEOREM 1.9. $\sup_{n \in \mathbb{N}} \mathbb{E}(\frac{n^2}{C[0,n]}) < \infty.$

1.4. *Outline.* Section 2 puts forward a structural representation of the invasion percolation cluster in terms of independent bond percolation, and gives the proof of Theorems 1.1 and 1.2. This structural representation plays a key role throughout the paper. Section 3 analyzes the process of forward maximal weights along the backbone and provides a scaling limit for this process in terms of the Poisson lower envelope process defined in (1.19). The multivariate Laplace transform of the latter is computed explicitly. Section 4 gives the proof of Theorem 1.3 and Corollary 1.4, based on the results in Section 3. Sections 5–8 give the proofs of Theorems 1.5, 1.6, 1.7 and 1.9, respectively. Section 9 provides a quick indication of proofs of the claims made in Section 1.2.

**2. Structural representation and local behavior.** In Section 2.1, we show that the IPC can be viewed as a random infinite backbone with subcritical percolation clusters emerging in all directions. The parameters of these subcritical clusters depend on the height of the vertex on the backbone from which they emerge, and tend to $p_c$ as this height tends to infinity. Theorem 1.1 follows immediately. In Section 2.2, we prove Theorem 1.2.

2.1. *Structural representation and proof of Theorem* 1.1.

2.1.1. *The structural representation.* As noted at the beginning of Section 1.3, the backbone is a.s. unique. Let $B_l$, $l \in \mathbb{N}$, denote the weights of its successive edges, and define

$$(2.1) \qquad W_k = \max_{l > k} B_l, \qquad k \in \mathbb{N}_0.$$

To see that the maximum in (2.1) is achieved, we first note that for each $k \in \mathbb{N}_0$ there must a.s. be an $l > k$ with $B_l > p_c$, since supercritical edges must be invaded to create an infinite cluster. On the other hand, we showed in Section 1.1 that for each $p > p_c$ there are at most finitely many edges invaded with weight above $p$. Thus the maximum in (2.1) is achieved, and $W_k > p_c$ a.s. In particular, $W_0$ is the weight of the heaviest edge on the backbone. Hence, it is also the weight of the heaviest edge ever invaded, since the existence of the infinite backbone path implies that no weight heavier than $W_0$ need ever be accepted.

The $W$-process is at the heart of our analysis, and we will study it in detail in Section 3. In particular, in a sense to be made precise in Proposition 3.3, we will see that

$$(2.2) \qquad W_k \sim p_c\left(1 + \frac{1}{k}Z\right) \qquad \text{as } k \to \infty$$

with $Z$ an exponential random variable with mean 1. This shows the slow rate of decay of $W_k$ toward the critical value.

The key observation behind the scaling results in Section 1.3 is the following structural representation of $C$ in terms of independent bond percolation.

PROPOSITION 2.1. *Under $\mathbb{P}$, $C$ can be viewed as consisting of:*

(1) *a single uniformly random infinite backbone;*
(2) *for all $k \in \mathbb{N}_0$, emerging from the $k$th vertex along the backbone, in all directions away from the backbone, an independent supercritical percolation cluster with parameter $W_k$ conditioned to stay finite.*

PROOF. By symmetry, all possible backbones are equally likely. We condition on the backbone, abbreviated BB. Conditional on $W = (W_k)_{k \in \mathbb{N}_0}$, the following is true for every $x \in \mathcal{T}_\sigma : x \in C$ if and only if every edge on the path between $x_{\text{BB}}$ and $x$ carries a weight below $W_k$, with $x_{\text{BB}}$ the vertex where the path downward from $x$ hits BB and $k = \|x_{\text{BB}}\|$. Indeed, if one of the edges in the path has weight above $W_k$, then this edge cannot be invaded, because the entire infinite BB is invaded first. Conversely, if all edges in the path have weight below $W_k$, then $x$ will be invaded before the edge on BB with weight $W_k$ is. In other words, the event $\{\text{BB} = bb, W = w\}$ is the same as the event that for all $k \in \mathbb{N}_0$ there is no percolation below level $W_k$ in each of the branches off BB at height $k$, and the forward maximal weights along $bb$ are equal to $w$. This proves the claim. $\square$

2.1.2. *The functions $\theta$ and $\zeta$.* For independent bond percolation on $\mathcal{T}_\sigma$ with parameter $p$, let $\theta(p)$ denote the probability that $o$ is in an infinite cluster, and let $\zeta(p)$ denote the probability that the cluster along a particular branch from $o$ is finite. Then we have the relations

$$(2.3) \qquad \theta(p) = 1 - \zeta(p)^\sigma, \qquad \zeta(p) = 1 - p\theta(p).$$

The critical probability is $p_c = 1/\sigma$, and $\theta(p_c) = 0$, $\zeta(p_c) = 1$.

For future reference, we note the following elementary facts. Differentiation of (2.3) gives

$$(2.4) \qquad \theta'(p) = \sigma \zeta(p)^{\sigma-1}[-\zeta'(p)], \qquad \zeta'(p) = -\theta(p) - p\theta'(p),$$

from which we see that

$$-\zeta'(p) = \frac{\theta(p)}{1 - p\sigma\zeta(p)^{\sigma-1}}. \tag{2.5}$$

The right-hand side gives $\frac{0}{0}$ for $p = p_c$. Using l'Hôpital's rule and the first equality of (2.4), we find that

$$-\zeta'(p_c) = \frac{\sigma[-\zeta'(p_c)]}{-\sigma + (\sigma - 1)[-\zeta'(p_c)]} \quad \text{and hence}$$
$$-\zeta'(p_c) = \frac{2\sigma}{\sigma - 1} = \frac{1}{\rho}, \tag{2.6}$$

where we recall the definition of $\rho$ in (1.3), and where derivatives at $p_c$ are interpreted as right-derivatives. From this, we obtain

$$\theta(p) \sim \frac{\sigma}{\rho}(p - p_c), \qquad 1 - \zeta(p) \sim \frac{1}{\rho}(p - p_c) \qquad \text{as } p \downarrow p_c. \tag{2.7}$$

In Section 3 we will need that $\zeta(p)$ is a convex function of $p \in [p_c, 1]$. This can be seen as follows. Since $\zeta$ is decreasing on $[p_c, 1]$ and maps this interval to $[0, 1]$, it is convex if and only if the inverse function $p = p(\zeta)$ is a convex function of $\zeta \in [0, 1]$. By (2.3), $p = F(\zeta)$ with $F(x) = \frac{1-x}{1-x^\sigma}$. Computation gives

$$F''(x) = \frac{\sigma x^{\sigma-2}}{(1 - x^\sigma)^3} G(x)$$
$$\text{with } G(x) = -(\sigma - 1)x^{\sigma+1} + (\sigma + 1)x^\sigma - (\sigma + 1)x + (\sigma - 1), \tag{2.8}$$

and hence it suffices to show that $G(x)$ is positive on $[0, 1]$. However, $G(1) = 0$, and

$$G'(x) = -(\sigma + 1)[-\sigma x^{\sigma-1} + (\sigma - 1)x^\sigma + 1] \tag{2.9}$$

is negative by the arithmetic-geometric mean inequality $(1 - \alpha)x_1 + \alpha x_2 \geq (x_1^{1-\alpha} x_2^\alpha)^{1/\alpha}$ with $\alpha = 1/\sigma$, $x_1 = x^\sigma$ and $x_2 = 1$.

For the special case $\sigma = 2$, (2.3) solves to give

$$\theta(p) = 0 \vee \frac{2p - 1}{p^2}, \qquad \zeta(p) = 1 \wedge \frac{1 - p}{p}. \tag{2.10}$$

2.1.3. *Duality and proof of Theorem* 1.1. The following duality is important in view of Proposition 2.1. Although this duality is standard in the theory of branching processes, we sketch the proof for completeness.

LEMMA 2.2. *On $\mathcal{T}_\sigma$, a supercritical percolation cluster with parameter $p > p_c$ conditioned to stay finite has the same law as a subcritical cluster with dual parameter*

$$\widehat{p} = \widehat{p}(p) = p\zeta(p)^{\sigma-1} < p_c. \tag{2.11}$$

*Moreover, $\widehat{p}(p_c) = p_c$, $\widehat{p}(1) = 0$, $\frac{d}{dp}\widehat{p}(p) < 0$ on $(p_c, 1)$, and*

$$p_c - \widehat{p}(p) \sim p - p_c \qquad \text{as } p \downarrow p_c. \tag{2.12}$$

For the special case $\sigma = 2$, (2.10) and (2.11) imply that the duality relation takes the simple form $\widehat{p} = 1 - p$.

PROOF OF LEMMA 2.2. Let $v$ be a vertex in $\mathcal{T}_\sigma$ and let $C(v)$ denote the forward cluster of $v$ for independent bond percolation with parameter $p$. Let $\mathcal{U}$ be any finite subtree of $\mathcal{T}_\sigma$, say with $m$ edges, and hence with $(\sigma - 1)m + \sigma$ boundary edges. Then

$$\begin{aligned}
\mathbb{P}_p(\mathcal{U} \subset C(v) \mid |C(v)| < \infty) &= \frac{\mathbb{P}_p(\mathcal{U} \subset C(v), |C(v)| < \infty)}{\mathbb{P}_p(|C(v)| < \infty)} \\
&= \frac{p^m \zeta(p)^{(\sigma-1)m+\sigma}}{\zeta(p)^\sigma},
\end{aligned} \tag{2.13}$$

the numerator being the probability that the edges of $\mathcal{U}$ are open and there is no percolation from any its vertices. Let

$$\widehat{p} = p\zeta(p)^{\sigma-1}. \tag{2.14}$$

Then the right-hand side of (2.13) equals $\widehat{p}^m = \mathbb{P}_{\widehat{p}}(\mathcal{U} \subset C(v))$. Since $\mathcal{U}$ is arbitrary, this proves the first claim.

Since the $\widehat{p}$ percolation clusters are a.s. finite we find $\widehat{p} \leq p_c$. Since $\zeta(p_c) = 1$ and $\zeta(1) = 0$, (2.14) implies that $\widehat{p}(p_c) = p_c$ and $\widehat{p}(1) = 0$. Direct computation gives $\frac{d}{dp}\widehat{p}(p) = \zeta(p)^{\sigma-1} + p(\sigma-1)\zeta(p)^{\sigma-2}\zeta'(p)$, which is negative if and only if $-\zeta'(p) > \zeta(p)/p(\sigma-1)$. By using (2.5) and (2.3), we see that the latter inequality holds if and only if $p\sigma > 1$, which is the same as $p > p_c$. Finally, we use the above formula for the derivative of $\widehat{p}(p)$, together with (2.6), to see that $\frac{d}{dp}\widehat{p}(p_c) = -1$ and hence

$$p_c - \widehat{p}(p) \sim p - p_c, \tag{2.15}$$

which is (2.12). □

Since a.s. $W_k > p_c$ for all $k \in \mathbb{N}_0$, we have $\widehat{W}_k < p_c$ for all $k \in \mathbb{N}_0$. Combining Proposition 2.1 and Lemma 2.2, we conclude that $C$ can be regarded as *a uniformly random infinite backbone with independent subcritical branches* with parameter $\widehat{W}_k$ emerging from the backbone vertex at height $k$ in all directions away from the backbone.

We are now in a position to better understand the difference between the IPC and the IIC. For the IIC, the branches emerging from the backbone are all critical percolation clusters. For the IPC, the branches are subcritical, and become increasingly close to critical as they branch off higher. Thus, low branches tend to be smaller than high branches. Conditional on $x \in C$, it is more likely for $x$ to be in a larger rather than a smaller branch, consistent with the observation in Section 1.3.2 that the backbone is more likely to branch off the path from $o$ to $x$ higher rather than lower.

The fact that the IPC is on average thinner than the IIC, as was observed in Section 1.3.3, is obvious from the fact that the subcritical branches of the IPC are smaller than the critical branches of the IIC. Moreover, the fact that there is randomness in the weights $\widehat{W}_k$ that determine the percolation parameters for the branches is consistent with the observation in Section 1.3.3 that the IPC has relatively larger fluctuations than the IIC.

PROOF OF THEOREM 1.1. It was noted in Section 1.2 that the IIC on $\mathcal{T}_\sigma$ can be viewed as consisting of a uniformly random infinite backbone with independent *critical* branches. In view of this observation, the statement made in Theorem 1.1 is an immediate consequence of Proposition 2.1 and Lemma 2.2. □

2.2. *Local behavior.*

PROOF OF THEOREM 1.2. The main idea in the proof is that a vertex $x \in C$ is unlikely to be very close to the backbone. On the other hand, the branch off the backbone containing $x$ is unlikely to branch close to $o$, and so it is close to critical percolation.

Fix a cylinder event $E$ on $\mathcal{T}_\sigma^*$. Let $k = k_E$ denote the maximal distance from $o$ to a vertex in a bond upon which $E$ depends. Fix $x \in \mathcal{T}_\sigma$. Let $M = M(x)$ denote the height of the highest vertex in the backbone on the path in $\mathcal{T}_\sigma$ from $o$ to $x$. As before, we write $W_M$ for the forward maximal weight above this vertex at height $M$ on the backbone. For $\varepsilon > 0$, let

$$
(2.16) \quad A^x = \{M \geq \|x\| - k\}, \qquad B^{x,\varepsilon} = \{W_M \geq p_c + \varepsilon\}
$$
$$
G^{x,\varepsilon} = (A^x \cup B^{x,\varepsilon})^c.
$$

It follows from (1.15) [although we have not yet proved (1.15), we will *not* use circular reasoning] that

$$
(2.17) \quad \lim_{\|x\| \to \infty} \mathbb{P}(A^x \mid x \in C) = 0 \qquad \forall \varepsilon > 0.
$$

We will prove that also

$$\lim_{\|x\|\to\infty} \mathbb{P}(B^{x,\varepsilon} \mid x \in C) = 0 \qquad \forall \varepsilon > 0, \tag{2.18}$$

implying

$$\lim_{\|x\|\to\infty} \mathbb{P}(G^{x,\varepsilon} \mid x \in C) = 1 \qquad \forall \varepsilon > 0. \tag{2.19}$$

To prove (2.18), we put $\|x\| = n$ and write

$$\mathbb{P}(B^{x,\varepsilon} \mid x \in C) = \sum_{m=0}^{n} \frac{\mathbb{P}(x \in C, M = m, B^{x,\varepsilon})}{\mathbb{P}(x \in C)}. \tag{2.20}$$

By (1.15), the denominator is at least $cn\sigma^{-n}$ for some $c > 0$. By Proposition 2.1 and Lemma 2.2, the numerator is at most $\sigma^{-m}[\widehat{p}(\varepsilon)]^{n-m}\mathbb{P}(W_m \geq p_c + \varepsilon)$ with $\widehat{p}(\varepsilon)$ the dual of $p_c + \varepsilon$ (we used the fact that $W_m \geq p$ implies $\widehat{W}_m \leq \widehat{p}$ for all $p > p_c$). Since $\widehat{p}(\varepsilon) \leq p_c = \sigma^{-1}$, we thus have

$$\mathbb{P}(B^{x,\varepsilon} \mid x \in C) \leq \frac{1}{cn} \sum_{m=0}^{n} \mathbb{P}(W_m \geq p_c + \varepsilon). \tag{2.21}$$

From Lemma 3.2 in Section 3.1 we will see that $\mathbb{P}(W_m \geq p_c + \varepsilon) \leq \exp[-c(\varepsilon)m]$ for all $m \in \mathbb{N}$ for some $c(\varepsilon) > 0$. Hence the sum in (2.21) is bounded in $n$ for fixed $\varepsilon$. This proves (2.18).

For each $\varepsilon > 0$, we have

$$\begin{aligned}|\mathbb{P}(\tau_x E \mid x \in C) - \mathbb{P}^*_\infty(E)| &\leq |\mathbb{P}(\tau_x E \mid x \in C) - \mathbb{P}(\tau_x E \mid x \in C, G^{x,\varepsilon})| \\ &\quad + |\mathbb{P}(\tau_x E \mid x \in C, G^{x,\varepsilon}) - \mathbb{P}^*_\infty(E)|.\end{aligned} \tag{2.22}$$

In view of (2.19), the first term on the right-hand side goes to zero as $\|x\| \to \infty$ for $\varepsilon > 0$ fixed, so it suffices to prove that

$$\lim_{\varepsilon \downarrow 0} \sup_{x \in \mathcal{T}^*_\sigma} |\mathbb{P}(\tau_x E \mid x \in C, G^{x,\varepsilon}) - \mathbb{P}^*_\infty(E)| = 0. \tag{2.23}$$

Now, on the event $\{x \in C\} \cap G^{x,\varepsilon}$, we have $\|x\| - k > M$, so that the event $\tau_x E$ depends only on bonds within a branch leaving the backbone at height $M$, and $W_M \in [p_c, p_c + \varepsilon)$, so that this branch is as close as desired to a critical tree when $\varepsilon$ is sufficiently small. Therefore, in the limit as $\varepsilon \downarrow 0$, $\mathbb{P}(\tau_x E \mid x \in C, G^{x,\varepsilon})$ approaches the probability of $E$ under the IIC rooted at $x$ and with a particular initial backbone segment of length $\|x\| - M$. The rate of convergence depends on the number of bonds upon which $E$ depends, but is uniform in $x$. However, by our hypothesis that $E$ is invariant under the automorphisms of $\mathcal{T}^*_\sigma$, $E$ has the same probability under the law $\mathbb{P}^*_\infty$ conditional on any choice of the initial backbone segment. This proves (2.23). $\square$

**3. Analysis of the backbone forward maximum process.** In this section, we prove that the backbone forward maximum process $W = (W_k)_{k \in \mathbb{N}_0}$ converges, after rescaling, to the Poisson lower envelope process $L = (L(t))_{t>0}$. In Section 3.1, we analyze $W$ as a Markov chain. In Section 3.2, we prove the convergence to $L$. Finally, in Section 3.3, we compute the multivariate Laplace transform of $L$.

3.1. *The Markov representation.*

PROPOSITION 3.1. *$W = (W_k)_{k \in \mathbb{N}_0}$ is a decreasing Markov chain taking values in $(p_c, 1)$ with initial distribution $\mathbb{P}(W_0 \leq u) = \theta(u)$ and transition probabilities*

$$\mathbb{P}(W_{k+1} = W_k \mid W_k = u) = 1 - R(u)\theta(u),$$
(3.1)
$$\mathbb{P}(W_{k+1} \in dv \mid W_k = u) = R(u)\theta'(v)\, dv,$$

*for $p_c < v < u < 1$, where $R(u) = \dfrac{1}{-\zeta'(u)}$.*

PROOF. The event $\{W_0 \leq u\}$ is the event that there is percolation at level $u$ on the tree, and hence has probability $\theta(u)$.

Denote by $\vec{W}_{<k}$ the vector $(W_j)_{0 \leq j < k}$. Clearly the process does not depend on which particular path forms the backbone, so we may fix the first $k$ edges of the backbone. Fix a vector $\vec{w}$ and $v \leq u \leq w_{k-1}$, and consider the conditional probability $\mathbb{P}(W_{k+1} \in dv \mid W_k = u, \vec{W}_{<k} = \vec{w})$. This is defined in terms of the conditional expectation

(3.2) $$\mathbb{E}[I(W_{k+1} \in dv) \mid W_k, \vec{W}_{<k}]$$

by setting $W_k = u$ and $\vec{W}_{<k} = \vec{w}$. We let $\vec{B}_{<k}$ denote the backbone weights below height $k$, and note that the above conditional expectation is equal to

(3.3) $$\mathbb{E}[\mathbb{E}[I(W_{k+1} \in dv) \mid W_k, \vec{B}_{<k}] \mid W_k, \vec{W}_{<k}],$$

since the pair $W_k, \vec{B}_{<k}$ specifies more information than the pair $W_k, \vec{W}_{<k}$. However, it is clear that

(3.4) $$\mathbb{E}[I(W_{k+1} \in dv) \mid W_k, \vec{B}_{<k}] = \mathbb{E}[I(W_{k+1} \in dv) \mid W_k],$$

since given $W_k$ the values of $\vec{B}_{<k}$ cannot affect $W_{k+1}$. Thus (3.2) is equal to

(3.5) $$\mathbb{E}[\mathbb{E}[I(W_{k+1} \in dv) \mid W_k] \mid W_k, \vec{W}_{<k}] = \mathbb{E}[I(W_{k+1} \in dv) \mid W_k].$$

This shows that $W$ is a Markov process.

To evaluate the transition probabilities we may consider only the case $k = 0$. We have already seen that

(3.6) $$\mathbb{P}(W_0 \in du) = \theta'(u)\, du.$$

For $v < u$, to have both $W_0 \in du$ and $W_1 \in dv$ there must also be an edge $e$ from the root such that:

1. The threshold for percolation above $e$ is in $dv$.
2. The weight of edge $e$ is $w_e \in du$.
3. There is no percolation at level $u$ in the other branches emerging from the root.

With $\sigma$ choices for $e$ we get

(3.7) $$\mathbb{P}(W_1 \in dv, W_0 \in du) = \sigma\theta'(v)\,dv\,du\zeta^{\sigma-1}(u).$$

Combining (3.6), (3.7) and using (2.4) we get

(3.8) $$\mathbb{P}(W_1 \in dv | W_0 = u) = \frac{\sigma\zeta^{\sigma-1}(u)}{\theta'(u)}\theta'(v)\,dv = R(u)\theta'(v)\,dv.$$

Finally, integrating over $v \in (p_c, u)$ we find

(3.9) $$\mathbb{P}(W_1 < W_0 \mid W_0 = u) = R(u)\theta(u),$$

and (3.1) follows from (3.8)–(3.9). $\square$

Note the separation in $u$ and $v$ in (3.1). The convexity of $\zeta$ (see Section 2.1.2) implies that $R$ is increasing and so, together with (2.6), yields

(3.10) $$R(u) \geq R(p_c) = \rho, u \in [p_c, 1].$$

For the special case $\sigma = 2$, (2.10) gives $R(u) = u^2$.

We have established already that $W_k > p_c$ for all $k \in \mathbb{N}_0$. The following large deviation estimate, which we applied in Section 2.2, shows that $W_k \downarrow p_c$ as $k \to \infty$, $\mathbb{P}$-a.s.

LEMMA 3.2. *For every $\delta > 0$ there is a $c(\delta) > 0$, satisfying $c(\delta) \sim \delta$ as $\delta \downarrow 0$, such that*

(3.11) $$\mathbb{P}\left(W_k \geq \frac{1}{\sigma}(1+\delta)\right) \leq e^{-c(\delta)k} \quad \text{and}$$
$$\mathbb{P}\left(\widehat{W}_k \leq \frac{1}{\sigma}(1-\delta)\right) \leq e^{-c(\delta)k} \qquad \forall k \in \mathbb{N}_0.$$

PROOF. We first claim that

(3.12) $$\mathbb{P}(W_k \geq p) \leq [1 - \rho\theta(p)]^k, \qquad k \in \mathbb{N}_0.$$

Indeed, (3.1) tells us that, for every $l \in \mathbb{N}_0$, given $W_l = u$, the probability that $W_{l+1} < p$ is $R(u)\theta(p)$. Hence, by (3.10), at each step $W$ has probability at least $\rho\theta(p)$ to jump below $p$, which implies (3.12). By (2.7), we have $\theta(\frac{1}{\sigma}(1+\delta)) \sim \frac{\delta}{\rho}$ as $\delta \downarrow 0$, and so we get the first part of (3.11). The second part follows from the first via Lemma 2.2. □

From (3.1), we have the following recursive representation for the $W$-process. Let $(X_k)_{k \in \mathbb{N}_0}$ be i.i.d. random variables with cumulative distribution function $\mathbb{P}(X_1 \leq u) = \theta(u)$, $u \in [0,1]$. Then $W_0 = X_0$ and, for $k \in \mathbb{N}_0$,

$$(3.13) \quad W_{k+1} = \begin{cases} W_k, & \text{with probability } 1 - R(W_k), \\ W_k \wedge X_{k+1}, & \text{with probability } R(W_k). \end{cases}$$

To prepare the ground for Proposition 3.3 below, let

$$(3.14) \quad Y_k = \rho \theta(W_k), \qquad k \in \mathbb{N}_0.$$

Note that $Y_k \downarrow 0$ as $k \to \infty$, $\mathbb{P}$-a.s., by Lemma 3.2. Let $(U_k)_{k \in \mathbb{N}_0}$ be i.i.d. uniform random variables on $[0,1]$. Then it follows from (3.13) that $Y = (Y_k)_{k \in \mathbb{N}_0}$ is a Markov chain with initial value $Y_0 = \rho U_0$ and recursive representation

$$(3.15) \quad Y_{k+1} = \begin{cases} Y_k, & \text{with probability } 1 - q(Y_k), \\ Y_k U_{k+1}, & \text{with probability } q(Y_k), \end{cases}$$

where

$$(3.16) \quad q(y) = \frac{y}{\rho} R\left(\theta^{-1}\left(\frac{y}{\rho}\right)\right)$$

with $\theta^{-1}$ the inverse of the function $\theta$. It then follows from (3.10) that

$$(3.17) \quad q(y) \geq y \quad \text{for } y \in [0, \rho] \quad \text{and} \quad q(y) \sim y \quad \text{as } y \downarrow 0.$$

This will be important in Section 3.2.

3.2. *Convergence of the forward maximum process to the Poisson lower envelope process.* The key to our analysis is the following proposition, which shows that the Poisson lower envelope process $L$ in (1.19) is the scaling limit of the backbone forward maximum process $W$ in (2.1). In particular, by taking $t = 1$ in (3.18) and using the fact that $L(1)$ is an exponential random variable with mean 1, we get the claim made in (2.2). We write $\stackrel{*}{\Longrightarrow}$ to denote convergence in distribution in the space of cadlag paths endowed with the Skorohod topology (see Billingsley [3], Section 14).

PROPOSITION 3.3. *For any $\varepsilon > 0$,*

$$(3.18) \quad (k[\sigma W_{\lceil kt \rceil} - 1])_{t \geq \varepsilon} \stackrel{*}{\Longrightarrow} (L(t))_{t \geq \varepsilon} \qquad \text{as } k \to \infty.$$

PROOF. The proof is based on the representation (3.15).

Let $N = (N(t))_{t \geq 0}$ denote the Poisson process on $[0, \infty)$ that increases at rate 1. Define

$$\widetilde{Y}(t) = Y_{N(t)}, \qquad t \geq 0. \tag{3.19}$$

Then $\widetilde{Y} = (\widetilde{Y}(t))_{t \geq 0}$ is the continuous-time Markov process with initial value $Y_0$ that from height $z$ jumps down to height $zU[0,1]$ at exponential rate $q(z)$. The $L$-process defined in (1.19) is the continuous-time Markov process that from height $z$ jumps down to height $zU[0,1]$ at exponential rate $z$. Below we will first use (3.17) to show that, for any $\varepsilon > 0$,

$$(k\widetilde{Y}(kt))_{t \geq \varepsilon} \overset{*}{\Longrightarrow} (L(t))_{t \geq \varepsilon}. \tag{3.20}$$

After that we will use the law of large numbers for $N$, namely $\lim_{k \to \infty} N(kt)/kt = 1$ a.s., to show that, for any $\varepsilon > 0$,

$$(kY_{\lceil kt \rceil})_{t \geq \varepsilon} \overset{*}{\Longrightarrow} (L(t))_{t \geq \varepsilon}. \tag{3.21}$$

Once we have (3.21), the proof is complete because

$$Y_{\lceil kt \rceil} \sim \sigma W_{\lceil kt \rceil} - 1 \qquad \text{as } k \to \infty \text{ uniformly in } t \geq \varepsilon, \tag{3.22}$$

as is immediate from (2.7) and (3.14), and the fact that the $Y$-process converges to 0, $\mathbb{P}$-a.s.

*Proof of* (3.20): The proof uses a *perturbative coupling argument*, relying on the fact that $q(z) \geq z$ for $z > 0$, while for every $\delta > 0$ there exists a $z_0 = z_0(\delta) > 0$ such that $q(z) \leq (1 + \delta)z$ for all $z \in (0, z_0]$.

*Upper bound*: For $y_0 > 0$, let $L_{y_0} = (L_{y_0}(t))_{t > 0}$ be the restriction to $(0, \infty) \times (0, y_0]$ of the lower envelope process associated with the Poisson process $\mathcal{P}$ (recall Figure 4), that is,

$$\begin{aligned} L_{y_0}(t) &= y_0 \wedge L(t) \\ &= y_0 \wedge \min\{y > 0 : (x, y) \in \mathcal{P} \text{ for some } x \leq t\}, \qquad t > 0. \end{aligned} \tag{3.23}$$

From height $z \leq y_0$, $L_{y_0}$ jumps down at exponential rate $z$. Therefore, conditional on $\widetilde{Y}_0 = y_0$, we can couple $\widetilde{Y}$ and $L_{y_0}$ such that

$$\widetilde{Y}(t) \leq L_{y_0}(t) \qquad \forall t > 0. \tag{3.24}$$

Indeed, to achieve the coupling we use the same uniform random variables for the jumps downward in both processes (so that the same sequence of heights are visited), but after each jump we arrange that $\widetilde{Y}$ waits less time than $L_{y_0}$ for its next jump, which is possible because $q(z) \geq z$ for $z > 0$. Combining (3.23) and (3.24), we find that $\widetilde{Y}$ and $L$ can be coupled so that

$$\widetilde{Y}(t) \leq L_{y_0}(t) \leq L(t) \qquad \forall t > 0. \tag{3.25}$$

This is a stochastic upper bound valid for all times.

*Lower bound*: We can imitate the above coupling argument, except that, in order to properly exploit the inequality $q(z) \leq (1+\delta)z$ for $z \in (0, z_0]$, we need a Poisson process with intensity $1 + \delta$, which we denote by $\mathcal{P}^{1+\delta}$, and we start the coupling only after $\widetilde{Y}$ has dropped below height $z_0$.

For $y_0 \leq z_0$, let $L_{y_0}^{1+\delta}$ be the restriction to $(0, \infty) \times (0, y_0]$ of the lower envelope process $L^{1+\delta}$ associated with $\mathcal{P}^{1+\delta}$, that is,

$$
\begin{aligned}
L_{y_0}^{1+\delta}(t) &= y_0 \wedge L^{1+\delta}(t) \\
&= y_0 \wedge \min\{y > 0 : (x, y) \in \mathcal{P}^{1+\delta} \text{ for some } x \leq t\}, \qquad t > 0.
\end{aligned}
\tag{3.26}
$$

Let

$$
T_0 = \min\{t > 0 : \widetilde{Y}(t) \leq z_0\}.
\tag{3.27}
$$

Then, conditional on $\widetilde{Y}(T_0) = y_0$, we can couple $\widetilde{Y}$ and $L_{y_0}^{1+\delta}$ in an analogous fashion to the coupling in the upper bound, such that

$$
\widetilde{Y}(t) \geq L_{y_0}^{1+\delta}(t) \qquad \forall t \geq T_0.
\tag{3.28}
$$

Next, let

$$
T_1 = \min\{x \geq T_0 : (x, y) \in \mathcal{P}^{1+\delta} \text{ for some } y < L_{y_0}^{1+\delta}(T_0)\}.
\tag{3.29}
$$

In words, $T_1$ is the first time after $T_0$ that $L_{y_0}^{1+\delta}(t)$ jumps. By construction,

$$
L_{y_0}^{1+\delta}(t) = L^{1+\delta}(t) \qquad \forall t \geq T_1.
\tag{3.30}
$$

Combining (3.28) and (3.30), we find that $\widetilde{Y}$ and $L^{1+\delta}$ can be coupled so that

$$
\widetilde{Y}(t) \geq L^{1+\delta}(t) \qquad \forall t \geq T_1.
\tag{3.31}
$$

This is a stochastic upper bound valid for large times, provided that $T_1 = T_1(T_0) < \infty$ a.s. For this to be true it suffices that $T_0 < \infty$ a.s. The latter is evidently true, because $q$ is bounded away from 0 outside any neighborhood of $z = 0$, implying that $\widetilde{Y}$ tends to 0 a.s.

*Sandwich*: For all $k \in \mathbb{N}$, $(kL(kt))_{t>0}$ has the same distribution as $(L(t))_{t>0}$, and $(L^{1+\delta}(t))_{t>0}$ has the same distribution as $(\frac{1}{1+\delta}L(t))_{t>0}$. Combined with (3.25) and (3.31), this implies that $\widetilde{Y}$ and $L$ can be coupled so that

$$
\frac{1}{1+\delta}L(t) \leq k\widetilde{Y}(kt) \leq L(t) \qquad \forall k \geq \widetilde{K} \text{ uniformly in } t \geq \varepsilon,
\tag{3.32}
$$

where $\widetilde{K} = \widetilde{K}(\delta, \varepsilon)$ is some finite random variable. Now let $\delta \downarrow 0$, to get the claim in (3.20).

*Proof of* (3.21): Fix $\gamma > 0$. Let $\mathcal{N}$ denote the law of $N$. By the strong law of large numbers, we have

(3.33)
$$N(kt) \in [(1-\gamma)kt, (1+\gamma)kt] \qquad \forall k \geq K \text{ uniformly in } t \geq \varepsilon, \mathcal{N}\text{-a.s.},$$

where $K = K(\gamma, \varepsilon)$ is some finite random variable. Because $\widetilde{Y}$ is a decreasing process, it follows from (3.19) and (3.33) that

(3.34)
$$kY_{\lceil kt \rceil} \in [k\widetilde{Y}((1+\gamma)kt), k\widetilde{Y}((1-\gamma)kt)]$$
$$\forall k \geq K \text{ uniformly in } t \geq \varepsilon, \mathbb{P} \times \mathcal{N}\text{-a.s.}$$

Combining (3.32) and (3.34), we find that there is a $K' = K'(\gamma, \delta, \varepsilon)$ such that $Y$ and $L$ can be coupled so that

(3.35) $\quad \dfrac{1}{1+\delta}\dfrac{1}{1+\gamma}L(t) \leq kY_{\lceil kt \rceil} \leq \dfrac{1}{1-\gamma}L(t) \qquad \forall k \geq K'$ uniformly in $t \geq \varepsilon$.

Now let $\delta, \gamma \downarrow 0$, to get the claim in (3.21). □

COROLLARY 3.4. *For any $\varepsilon > 0$,*

(3.36)
$$(k[1 - \sigma\widehat{W}_{\lceil kt \rceil}])_{t \geq \varepsilon} \stackrel{*}{\Longrightarrow} (L(t))_{t \geq \varepsilon} \qquad \text{as } k \to \infty.$$

PROOF. By Lemma 3.2, $W_k \downarrow p_c$ as $k \to \infty$, $\mathbb{P}$-a.s., so (3.36) is immediate from (2.12) and (3.18). □

3.3. *Multivariate Laplace transform of the Poisson lower envelope.* Recall the definition of the $L$-process in (1.19). The following lemma gives its multivariate Laplace transform.

LEMMA 3.5. *For any $n \in \mathbb{N}$, $\tau_1, \ldots, \tau_n \geq 0$ and $0 \leq t_1 < \cdots < t_n$,*

(3.37)
$$\mathbb{E}\left(\exp\left[-\sum_{i=1}^n \tau_i L(t_i)\right]\right) = \prod_{i=1}^n \left(1 - \frac{\tau_i}{t_i + s_i}\right)$$

*with $s_i = \sum_{j=1}^i \tau_j$.*

PROOF. Let

(3.38)
$$\mathcal{I} = \{1 \leq i < n : L(t_{i+1}) < L(t_i)\}.$$

We split the contribution according to the outcome of $\mathcal{I}$. To that end, fix $0 \leq m \leq n-1$ and $A = \{a_1, \ldots, a_m\}$ with $1 \leq a_1 < \cdots < a_m \leq n-1$. Put $a_0 = 0$ and $a_{m+1} = n$. On the event $\{\mathcal{I} = A\}$, there are $u_1 > u_2 > \cdots > u_m > u_{m+1} > 0$ such that

(3.39) $\forall j = 1, \ldots, m+1$: $L(t_i) \in (u_j, u_j + du_j]$ for $a_{j-1} < i \leq a_j$.

In terms of the Poisson process $\mathcal{P}$, this is the same as the event

$$(3.40) \quad \forall j = 1, \ldots, m+1 : \begin{cases} \mathcal{P} \cap (t_{a_{j-1}}, t_{a_j}] \times (0, u_j] = \varnothing, \\ \mathcal{P} \cap (t_{a_{j-1}}, t_{a_{j-1}+1}] \times (u_j, u_j + du_j] \neq \varnothing, \end{cases}$$

where we put $t_0 = 0$. The latter event has probability

$$(3.41) \quad \prod_{j=1}^{m+1} e^{-u_j(t_{a_j} - t_{a_{j-1}})} (t_{a_{j-1}+1} - t_{a_{j-1}}) \, du_j.$$

Furthermore, on this event we have

$$(3.42) \quad \sum_{i=1}^{n} \tau_i L(t_i) = \sum_{j=1}^{m+1} u_j (s_{a_j} - s_{a_{j-1}}),$$

where we put $s_0 = 0$. Therefore, we obtain

$$(3.43) \quad \begin{aligned} & \mathbb{E}\left( \exp\left[ -\sum_{i=1}^{n} \tau_i L(t_i) \right] 1_{\{\mathcal{I}=A\}} \right) \\ & = \left( \prod_{j=1}^{m+1} \int_0^\infty du_j \right) 1_{\{u_1 > u_2 > \cdots > u_m > u_{m+1}\}} \\ & \quad \times \prod_{j=1}^{m+1} (t_{a_{j-1}+1} - t_{a_{j-1}}) e^{-u_j[(t_{a_j} - t_{a_{j-1}}) + (s_{a_j} - s_{a_{j-1}})]}. \end{aligned}$$

It is straightforward to perform the integrals in (3.43) in the order $j = 1, \ldots, m+1$, noting that the exponent telescopes, to get

$$(3.44) \quad \prod_{j=1}^{m+1} \frac{(t_{a_{j-1}+1} - t_{a_{j-1}})}{(t_{a_j} - t_{a_0}) + (s_{a_j} - s_{a_0})}.$$

Since $a_0 = 0$, $t_0 = 0$ and $s_0 = 0$, this gives the formula

$$(3.45) \quad \begin{aligned} & \mathbb{E}\left( \exp\left[ -\sum_{i=1}^{n} \tau_i L(t_i) \right] 1_{\{\mathcal{I}=A\}} \right) \\ & = \prod_{j=1}^{m+1} \frac{t_{a_{j-1}+1} - t_{a_{j-1}}}{t_{a_j} + s_{a_j}} = \frac{t_1}{t_n + s_n} \prod_{i \in A} \frac{t_{i+1} - t_i}{t_i + s_i}, \end{aligned}$$

with the empty product equal to 1. Finally, we sum over $A$ and use that

$$(3.46) \quad \sum_A \prod_{i \in A} \frac{t_{i+1} - t_i}{t_i + s_i} = \prod_{i=1}^{n-1} \left( 1 + \frac{t_{i+1} - t_i}{t_i + s_i} \right) = \prod_{i=1}^{n-1} \frac{t_{i+1} + s_i}{t_i + s_i},$$

to arrive at

$$\mathbb{E}\left(\exp\left[-\sum_{i=1}^{n}\tau_i L(t_i)\right]\right) = \prod_{i=1}^{n}\frac{t_i + s_{i-1}}{t_i + s_i}, \qquad (3.47)$$

which is the formula in (3.37). □

## 4. Proof of Theorem 1.3 and Corollary 1.4.

PROOF OF THEOREM 1.3. For fixed $\vec{x}$, and for $w, v \in \mathcal{N}(\vec{x})$ with $w \leq v$, let $N_w^v$ denote the number of edges in the connected component of $w$ in the subgraph of $\mathcal{S}(\vec{x})$ that is obtained by removing all the edges in the path from $o$ to $v$ (see Figure 2). Pick $y \in \partial \mathcal{S}(\vec{x})$, and let $v \in \mathcal{N}(\vec{x}) \setminus \{o\}$ be the first node above the parent of $y$, and let $0 < k \leq n_v$ be the distance from $v_-$ to the parent of $y$. Then the event $\{\vec{x} \in C\} \cap B_y$ amounts to the following:

(1) The backbone runs from $o$ to $v_-$, runs up a height $k$ along the segment between $v_-$ and $v$, and then moves to $y$;

(2) for all $w \in \mathcal{N}(\vec{x})$ with $w < v$, $N_w^v$ invaded edges are connected to the backbone at height $\|w\|$;

(3) $N_v^v + (n_v - k)$ invaded edges are connected to the backbone at height $\|v_-\| + k$.

Therefore

$$(4.1) \quad \begin{aligned}&\mathbb{P}(\vec{x} \in C, B_y \mid W) \\ &= \left(\frac{1}{\sigma}\right)^{\|v_-\|+k+1} \left\{\prod_{\substack{w \in \mathcal{N}(\vec{x}) \\ w < v}} [\widehat{W}_{\|w\|}]^{N_w^v}\right\} [\widehat{W}_{\|v_-\|+k}]^{N_v^v + (n_v - k)},\end{aligned}$$

where the three factors correspond to (1)–(3), and Proposition 2.1 and Lemma 2.2 are used to determine the probabilities of (2) and (3). Taking the average over $W$ and using the relation

$$\|v_-\| + k + \sum_{w \leq v} N_w^v + (n_v - k) = N, \qquad (4.2)$$

we obtain

$$(4.3) \quad \sigma^{N+1}\mathbb{P}(\vec{x} \in C, B_y) = \mathbb{E}\left(\left\{\prod_{\substack{w \in \mathcal{N}(\vec{x}) \\ w < v}} [\sigma\widehat{W}_{\|w\|}]^{N_w^v}\right\} [\sigma\widehat{W}_{\|v_-\|+k}]^{N_v^v + (n_v - k)}\right).$$

Since, by assumption, $\mathcal{S}(\vec{x})$ does not branch at $o$, we have $N_o^v = 0$, and so the factor with $w = o$ may be dropped.

We next apply Corollary 3.4 in combination with the scaling limit defined by (1.8)–(1.9). To that end, we define

$$
(4.4) \quad h_w(N) = \frac{\|w\|}{N}, \qquad n_w^v(N) = \frac{N_w^v}{N}, \qquad s(N) = \frac{k}{N},
$$
$$
t_v(N) = \frac{n_v}{N}, \qquad Z_h(N) = N[1 - \sigma \widehat{W}_{hN}],
$$

and $f_N(x) = (1 - x/N)^N 1_{[0,N]}(x)$, $x \in (0, \infty)$, and rewrite the right-hand side of (4.3) as

$$
(4.5) \quad \mathbb{E}\left(\left\{\prod_{\substack{w \in \mathcal{N}(\vec{x}) \\ o < w < v}} [f_N(Z_{h_w(N)}(N))]^{n_w^v(N)}\right\} \right.
$$
$$
\left. \times [f_N(Z_{h_{v_-}(N)+s(N)}(N))]^{n_v^v(N)+[t_v(N)-s(N)]}\right).
$$

Under the limit $\sharp \lim_{N \to \infty}$, there are $h_w$, $n_w^v$, $s$ and $t_v$ such that

$$(4.6) \quad h_w(N) \to h_w, \qquad n_w^v(N) \to n_w^v, \qquad s(N) \to s, \qquad t_v(N) \to t_v,$$

and, by Corollary 3.4,

$$
(4.7) \quad \begin{array}{c} ((Z_{h_w(N)}(N))_{o<w<v}, Z_{h_{v_-}(N)+s(N)}(N)) \\ \Longrightarrow ((L(h_w))_{o<w<v}, L(h_{v_-} + s)), \end{array}
$$

provided we assume that $s > 0$ when $v_- = o$. This last assumption (which will be removed below) is needed here because Corollary 3.4 only applies for positive scaling heights. Let $f(x) = \exp(-x)$, $x \in (0, \infty)$. Since $f_N$ converges to $f$ as $N \to \infty$ uniformly on $(0, \infty)$, and since $f$ is bounded and continuous on $(0, \infty)$, it follows from (4.5)–(4.7) that

$$
(4.8) \quad \begin{aligned} &\sharp \lim_{N \to \infty} \sigma^{N+1} \mathbb{P}(\vec{x} \in C, B_y) \\ &= \mathbb{E}\left(\left\{\prod_{\substack{w \in \mathcal{N} \\ o<w<v}} [f(L(h_w))]^{n_w^v}\right\} [f(L(h_{v_-} + s))]^{[n_v^v + (t_v - s)]}\right) \\ &= \mathbb{E}\left(\exp\left[-\sum_{\substack{w \in \mathcal{N} \\ o<w<v}} n_w^v L(h_w) - [n_v^v + (t_v - s)]L(h_{v_-} + s)\right]\right), \end{aligned}
$$

where $\mathcal{N}$ is the set of nodes in the scaled spanning tree $\mathcal{S}$ (defined above Theorem 1.3). Next, we use Lemma 3.5 to obtain

$$(4.9) \quad \sharp \lim_{N \to \infty} \sigma^{N+1} \mathbb{P}(\vec{x} \in C, B_y) = \left\{\prod_{\substack{w \in \mathcal{N} \\ o<w<v}} \left(1 - \frac{n_w^v}{m_w^v}\right)\right\}(1 - [n_v^v + (t_v - s)]),$$

where we recall the definition of $m_w^v$ in (1.10), and use the relations $h_w + \sum_{o<u\leq w} n_u^v = m_w^v$ and $h_{v_-} + s + \sum_{o<w\leq v_-} n_w^v + n_v^v + (t_v - s) = 1$ [by (4.2) with $N_o^v = 0$]. Finally, we note that $m_w^v - n_w^{\bar{v}} = t_w + m_{w_-}^v$ and $1 - (n_v^v + t_v) = m_{v_-}^v$, to obtain the formula in (1.11).

It is easy to remove the restriction that $s > 0$ when $v_- = o$. Indeed, the right-hand side of (4.3) is increasing in $k$, because $\widehat{W}_l$ is increasing in $l$. Therefore, we can include the case $s = 0$ via a monotone limit. $\square$

PROOF OF COROLLARY 1.4. In the limit as $N \to \infty$, the sum over $0 < k \leq n_v$ may be replaced by an integral over $s \in [0, t_v]$ for all $v \in \mathcal{N}\setminus\{o\}$, by using the monotonicity in $k$ noted above. $\square$

REMARK. If $\mathcal{S}(\vec{x})$ branches at $o$, then the factor $[\sigma\widehat{W}_0]^{N_o^v}$ in (4.3) is not 1. In fact, it tends to zero as $N \to \infty$, $\mathbb{P}$-a.s., because $\sigma\widehat{W}_0$ is a random variable on $(0, 1)$ while $N_o^v \sim n_o^v N \to \infty$ since $n_o^v > 0$. Thus the right-hand side of (4.3) tends to zero in this case.

**5. Cluster size at a given height.** In this section, we prove Theorem 1.5. The cluster size at height $n$ consists of the contributions at height $n$ from branches leaving the backbone at height $k$, for all $0 \leq k < n$, plus the single backbone vertex at height $n$. This leads us, in Section 5.1, to first analyze the Laplace transform of $C_o[m]$, which is the contribution to the cluster at height $m$ from a single branch from the root, in independent bond percolation with parameter $p$. Section 5.2 then uses this Laplace transform to provide the proof of Theorem 1.5, while Section 5.3 computes the first and second moment in the scaling limit.

5.1. *Laplace transform of $C_o[m]$.* For $m \in \mathbb{N}$, let $C_o[m]$ denote the number of vertices in the cluster of $o$ at height $m$, via a fixed branch from $o$, in independent bond percolation with parameter $p \leq p_c = 1/\sigma$. For $\tau \geq 0$, let

$$(5.1) \qquad f_m(p;\tau) = \mathbb{E}_p(e^{-\tau C_o[m]}).$$

By conditioning on the occupation status of the edge leaving the root, we see that $C_o[m+1]$ is 0 with probability $1-p$ and is the sum of $\sigma$ independent copies of $C_o[m]$ with probability $p$. Therefore $f_m$ obeys the recursion relation

$$(5.2) \quad f_{m+1}(p;\tau) = 1 - p + p[f_m(p;\tau)]^\sigma, \qquad f_1(p;\tau) = 1 - p + pe^{-\tau}.$$

We set $f_0(p;\tau) = e^{-\tau/\sigma}$, so that the recursion in (5.2) holds for $m = 0$ as well.

Let $g_m(p;\tau) = 1 - f_m(p;\tau)$ and $\rho = \frac{\sigma-1}{2\sigma}$. Our goal is to determine the asymptotic behavior of $g_m(p;\tau)$ for small $\tau$ and for $p$ near $p_c$. To emphasize the latter, we sometimes write $p = \frac{1}{\sigma}(1-\delta)$. However, we usually suppress the arguments $p$ and $\tau$. In terms of $g_m$, the recursion reads

$$(5.3) \qquad g_{m+1} = F(g_m), \qquad g_0 = 1 - e^{-\tau/\sigma},$$

with

$$(5.4) \qquad F(x) = p[1 - (1-x)^\sigma], \qquad x \in [0,1].$$

We will first show that the sequence $(g_m)_{m \in \mathbb{N}_0}$ is close to the sequence $(\widehat{g}_m)_{m \in \mathbb{N}_0}$ that satisfies the quadratic recursion

$$(5.5) \qquad \widehat{g}_{m+1} = \widehat{F}(\widehat{g}_m), \qquad \widehat{g}_0 = g_0,$$

where $\widehat{F}(x)$ is the second-order approximation of $F(x)$, namely,

$$(5.6) \qquad \widehat{F}(x) = p\sigma\left[x - \frac{\sigma-1}{2}x^2\right], \qquad x \in [0,1].$$

After that we will show that the sequence $(\widehat{g}_m)_{m \in \mathbb{N}_0}$ is close to the sequence $(\widetilde{g}(m))_{m \in \mathbb{N}_0}$, where $\widetilde{g}(t)$ satisfies the differential equation

$$(5.7) \qquad \widetilde{g}'(t) = \widehat{F}(\widetilde{g}(t)) - \widetilde{g}(t), \qquad \widetilde{g}(0) = g_0.$$

The differential equation (5.7) is easily solved, as follows. We abbreviate $q = \frac{1-\delta}{\delta}p\sigma$, where $p = \frac{1}{\sigma}(1-\delta)$, and rewrite (5.7) as

$$(5.8) \qquad \left(\frac{1}{\widetilde{g}} - \frac{q}{1+q\widetilde{g}}\right)\widetilde{g}' = -\delta.$$

This may be integrated to give

$$(5.9) \qquad \frac{\widetilde{g}(t)}{1+q\widetilde{g}(t)} = \frac{g_0}{1+qg_0}e^{-\delta t} \quad \text{or} \quad \widetilde{g}(t) = \frac{g_0 e^{-\delta t}}{1+qg_0[1-e^{-\delta t}]}.$$

The following lemma bounds the difference between $g_m$ and $\widetilde{g}(m)$.

LEMMA 5.1. *For $m \in \mathbb{N}_0$, $p \leq p_c$ and $\tau \geq 0$,*

$$(5.10) \qquad -\frac{1}{\sigma}\left(\delta\tau + \frac{1}{2}\tau^2\right) \leq g_m(p;\tau) - \widetilde{g}(m) \leq \frac{1}{6\sigma}m\tau^3,$$

*with $\widetilde{g}(m)$ given by (5.9).*

PROOF. We will prove the sandwiches

$$(5.11) \qquad 0 \leq g_m - \widehat{g}_m \leq \frac{1}{6\sigma}m\tau^3$$

and

$$(5.12) \qquad 0 \leq \widetilde{g}(m) - \widehat{g}_m \leq \frac{1}{\sigma}\left(\delta\tau + \frac{1}{2}\tau^2\right),$$

which together give the lemma.

We begin with (5.11). Since $0 \leq F'''(x) \leq p\sigma^3 \leq \sigma^2$, it follows from a third-order Taylor expansion that

$$(5.13) \qquad 0 \leq F(x) - \widehat{F}(x) \leq \tfrac{1}{6}\sigma^2 x^3, \qquad x \in [0,1].$$

Moreover, $F(0) = \widehat{F}(0) = 0$, $0 \leq F'(x) \leq 1$, and $\widehat{F}'(x) \leq 1$, so $0 \leq F(x) \leq x$ and $\widehat{F}(x) \leq x$ for all $x \in [0,1]$, and hence $g_m$, $\widehat{g}_m$ and $\widetilde{g}(t)$ are all decreasing. Write

$$(5.14) \qquad g_{m+1} - \widehat{g}_{m+1} = [F(g_m) - F(\widehat{g}_m)] + [F(\widehat{g}_m) - \widehat{F}(\widehat{g}_m)].$$

If $g_m \geq \widehat{g}_m$, then $F(g_m) - F(\widehat{g}_m) \geq 0$ by the monotonicity of $F$, while $F(\widehat{g}_m) - \widehat{F}(\widehat{g}_m) \geq 0$ because $F \geq \widehat{F}$, and so $g_{m+1} \geq \widehat{g}_{m+1}$. Since $g_0 = \widehat{g}_0$, it follows inductively that

$$(5.15) \qquad g_m \geq \widehat{g}_m \qquad \forall m \in \mathbb{N}_0.$$

Moreover, $F(g_m) - F(\widehat{g}_m) \leq g_m - \widehat{g}_m$ because $F' \leq 1$, while $F(\widehat{g}_m) - \widehat{F}(\widehat{g}_m) \leq \tfrac{1}{6}\sigma^2 \widehat{g}_m^3$ by (5.13). Therefore $g_{m+1} - \widehat{g}_{m+1} \leq g_m - \widehat{g}_m + \tfrac{1}{6}\sigma^2 \widehat{g}_m^3$. Since $\widehat{g}_m$ is decreasing and $\widehat{g}_0 = g_0 \leq \tau/\sigma$, this yields (5.11).

Next we prove (5.12). Define $\widetilde{F}(x) = h(1)$, where $h = h_x$ is the solution of

$$(5.16) \qquad h'(t) = \widehat{F}(h(t)) - h(t), \qquad h(0) = x.$$

According to (5.7), we have

$$(5.17) \qquad \widetilde{g}(m+1) = \widetilde{F}(\widetilde{g}(m)).$$

Since $\widehat{F}(x) \leq x$, the solution of (5.16) is decreasing, and therefore the function $\widetilde{F}$ is increasing. Now,

$$(5.18) \qquad h(1) - h(0) = \int_0^1 dt [\widehat{F}(h(t)) - h(t)],$$

and so, $h(t)$ and $\widehat{F}(x) - x$ both being decreasing, we have

$$(5.19) \qquad \widehat{F}(h(0)) - h(0) \leq h(1) - h(0) \leq \widehat{F}(h(1)) - h(1).$$

Since $h(1) = \widetilde{F}(h(0))$, the lower bound in (5.19) with $h(0) = x$ gives $\widetilde{F}(x) \geq \widehat{F}(x)$. With $\widetilde{g}(0) = \widehat{g}_0$, because of (5.5) and (5.17) and the fact that $\widetilde{F}$ is increasing, the latter inductively implies that

$$(5.20) \qquad \widetilde{g}(m) \geq \widehat{g}_m \qquad \forall m \in \mathbb{N}_0.$$

Using the upper bound in (5.19) with $h(0) = \widetilde{g}(k-1)$ and $h(1) = \widetilde{g}(k)$, and once more that $\widehat{F}(x) - x$ is decreasing in combination with (5.20), we get

$$(5.21) \quad \widetilde{g}(k) - \widetilde{g}(k-1) \leq \widehat{F}(\widetilde{g}(k)) - \widetilde{g}(k) \leq \widehat{F}(\widehat{g}_k) - \widehat{g}_k = \widehat{g}_{k+1} - \widehat{g}_k.$$

Summing (5.21) from $k=1$ to $k=m-1$, we obtain $\widetilde{g}(m-1)-\widetilde{g}(0) \leq \widehat{g}_m - \widehat{g}_1$. Since $\widetilde{g}(m) \leq \widetilde{g}(m-1)$, $\widetilde{g}(0) = g_0$ and $\widehat{g}_1 = \widehat{F}(\widehat{g}_0) = \widehat{F}(g_0)$, this yields the sandwich

$$(5.22) \qquad 0 \leq \widetilde{g}(m) - \widehat{g}_m \leq g_0 - \widehat{F}(g_0) \leq \frac{1}{\sigma}\left(\delta\tau + \frac{1}{2}\tau^2\right),$$

where the last inequality is for $p = \frac{1}{\sigma}(1-\delta)$ and uses $g_0 \leq \tau/\sigma$. This proves (5.12). $\square$

5.2. *Proof of Theorem* 1.5. Let $C_{k,j}[n]$ denote the contribution to the cluster size at height $n$ due to the $j$th of the $\sigma - 1$ branches emerging from height $k$ on the backbone. Then $C[n] = 1 + \sum_{k=0}^{n-1} \sum_{j=1}^{\sigma-1} C_{k,j}[n]$, with the additional 1 due to the backbone vertex at level $n$. We note that $C_{k,j}[n]$ is a random functional of $W_k$ and that, *conditional* on $W = (W_k)_{k \in \mathbb{N}_0}$, the different clusters are all independent. According to Proposition 2.1 and Lemma 2.2, each cluster emerging from the backbone at height $k$ is a subcritical cluster with parameter $p = \widehat{W}_k$. Thus,

$$(5.23) \qquad \mathbb{E}(e^{(-\tau/(\rho n))C_{k,j}[n]} \mid W) = f_{n-k}\left(\widehat{W}_k; \frac{\tau}{\rho n}\right),$$
$$k = 0, \ldots, n-1; j = 1, \ldots, \sigma-1; \tau \geq 0,$$

with $f_m(p;\tau)$ as defined in (5.1). Consequently,

$$\mathbb{E}(e^{-\tau\Gamma_n} \mid W) = \mathbb{E}(e^{(-\tau/(\rho n))C[n]} \mid W)$$
$$(5.24) \qquad = e^{(-\tau/(\rho n))} \left(\prod_{k=0}^{n-1} f_{n-k}\left(\widehat{W}_k; \frac{\tau}{\rho n}\right)\right)^{\sigma-1}$$
$$= e^{(-\tau/(\rho n))} \exp\left[(\sigma-1)\sum_{k=0}^{n-1} \log\left(1 - g_{n-k}\left(\widehat{W}_k; \frac{\tau}{\rho n}\right)\right)\right].$$

Note that, compared to Section 5.1, the argument $\tau$ has now become $\tau/\rho n$.

Fix $\tau \geq 0$. Fix $\varepsilon > 0$ and, for notational simplicity, assume that $\varepsilon n$ is integer. Since $g_{n-k} \leq g_0 \leq \tau/\sigma\rho n$, bounding the first $\varepsilon n + 1$ and the last $\varepsilon n - 1$ terms in the sum individually and using a linear approximation of the logarithm for the remaining terms, we get

$$(5.25) \qquad \mathbb{E}(e^{-\tau\Gamma_n} \mid W) = \exp\left[-S_n^\varepsilon(\tau, W) + O\left(\varepsilon\tau + \frac{\tau^2}{n}\right)\right]$$

with

$$(5.26) \qquad S_n^\varepsilon(\tau, W) = (\sigma-1) \sum_{\varepsilon n < k \leq n-\varepsilon n} g_{n-k}\left(\widehat{W}_k; \frac{\tau}{\rho n}\right).$$

Next, for $t \in [\varepsilon, 1-\varepsilon]$, put $k = \lceil nt \rceil$ and $Z_n^t = n[1 - \sigma \widehat{W}_{\lceil nt \rceil}]$, and define

$$(5.27) \quad G_n^t(x) = (\sigma-1)n g_{\lfloor n(1-t) \rfloor}\left(\frac{1}{\sigma}\left(1 - \frac{x}{n}\right); \frac{\tau}{\rho n}\right) 1_{[0,n]}(x), \qquad x \in (0, \infty).$$

Then we may write (5.26) as

$$(5.28) \quad S_n^\varepsilon(\tau, W) = \int_\varepsilon^{1-\varepsilon} dt\, G_n^t(Z_n^t).$$

Applying Lemma 5.1, after some arithmetic we find that for $x \in [0, n]$ and $n \to \infty$,

$$(5.29) \quad \begin{aligned} G_n^t(x) &= \left[1 + O\left(\frac{\tau}{n}\right)\right] 2\tau \frac{xe^{-(1-t)x}}{x + [1 + O((x+\tau)/n)]\tau[1 - e^{-(1-t)x}]} \\ &\quad + O\left(\frac{x\tau + \tau^2 + \tau^3}{n}\right). \end{aligned}$$

Put

$$(5.30) \quad G^t(x) = 2\tau \frac{xe^{-(1-t)x}}{x + \tau[1 - e^{-(1-t)x}]}, \qquad x \in (0, \infty).$$

Note that $G^t \leq 2\tau$ and that $\lim_{x \to \infty} G^t(x) = 0$ uniformly in $t \in [\varepsilon, 1-\varepsilon]$. Also, it follows from the fact that $f_m(p; \tau)$ is decreasing in $p$ that $G_n^t(x)$ is decreasing in $x$. It is clear from (5.29) that $G_n^t(x) \to G^t(x)$ uniformly in $x \in (0, \sqrt{n}]$ and $t \in (\varepsilon, 1-\varepsilon)$. For $x > \sqrt{n}$, the difference $|G_n^t(x) - G^t(x)|$ is bounded above by $G_n^t(\sqrt{n}) + G^t(x)$, which also goes to zero uniformly in $x > \sqrt{n}$ and $t \in (\varepsilon, 1-\varepsilon)$. Therefore, $G_n^t(x)$ converges to $G^t(x)$ as $n \to \infty$, uniformly in $x \in (0, \infty)$ and $t \in [\varepsilon, 1-\varepsilon]$.

Consequently,

$$(5.31) \quad \lim_{n \to \infty} \int_\varepsilon^{1-\varepsilon} dt[G_n^t(Z_n^t) - G^t(Z_n^t)] = 0, \qquad \mathbb{P}\text{-a.s.}$$

Moreover, we know from Corollary 3.4 that $(Z_n^t)_{t \in [\varepsilon, 1-\varepsilon]} \overset{*}{\Longrightarrow} (L(t))_{t \in [\varepsilon, 1-\varepsilon]}$ as $n \to \infty$. Since $(z(t))_{t \in [\varepsilon, 1-\varepsilon]} \mapsto \int_\varepsilon^{1-\varepsilon} dt\, G^t(z(t))$ is bounded and continuous in the Skorohod topology, it follows that

$$(5.32) \quad \int_\varepsilon^{1-\varepsilon} dt\, G^t(Z_n^t) \Longrightarrow \int_\varepsilon^{1-\varepsilon} dt\, G^t(L(t)) \qquad \text{as } n \to \infty.$$

Combining (5.28), (5.31) and (5.32), we obtain

$$(5.33) \quad S_n^\varepsilon(\tau, W) \Longrightarrow \int_\varepsilon^{1-\varepsilon} dt\, G^t(L(t)) \qquad \text{as } n \to \infty.$$

We substitute (5.33) into (5.25), and let $n \to \infty$ followed by $\varepsilon \downarrow 0$, to get

$$(5.34) \quad \lim_{n \to \infty} \mathbb{E}(e^{-\tau \Gamma_n}) = \mathbb{E}\left(\exp\left[-\int_0^1 dt\, G^t(L(t))\right]\right), \qquad \tau \geq 0.$$

The integral in the right-hand side equals $S(\tau, L)$ of (1.21), and this proves (1.20).

5.3. *First and second moment.* In this section we prove (1.22).
Differentiation of (1.21) gives

$$(5.35) \quad \mathbb{E}(\Gamma) = \mathbb{E}\left(\frac{\partial S}{\partial \tau}(0, L)\right) = \mathbb{E}\left(2\int_0^1 e^{-(1-t)L(t)}\,dt\right) = 2\int_0^1 t\,dt = 1,$$

where the third equality uses Lemma 3.5. Similarly,

$$(5.36) \quad \mathbb{E}(\Gamma^2) = I + II$$

with

$$(5.37) \quad \begin{aligned} I &= \mathbb{E}\left(\left[\frac{\partial S}{\partial \tau}(0, L)\right]^2\right) = \mathbb{E}\left(4\int_0^1 dt \int_0^1 ds\, e^{-(1-t)L(t)} e^{-(1-s)L(s)}\right) \\ &= \mathbb{E}\left(8\int_0^1 dt \int_t^1 ds\, e^{-(1-t)L(t)} e^{-(1-s)L(s)}\right) \\ &= 8\int_0^1 dt \int_t^1 ds\, t\frac{1-t+s}{2-t} \\ &= 4\int_0^1 t\left((2-t) - \frac{1}{2-t}\right) dt = \frac{20}{3} - 8\log 2 \end{aligned}$$

and

$$(5.38) \quad \begin{aligned} II &= \mathbb{E}\left(-\frac{\partial^2 S}{\partial \tau^2}(0, L)\right) = \mathbb{E}\left(4\int_0^1 \frac{dt}{L(t)} e^{-(1-t)L(t)}\left[1 - e^{-(1-t)L(t)}\right]\right) \\ &= \mathbb{E}\left(4\int_0^1 dt \int_t^1 ds\, e^{-(2-t-s)L(t)}\right) \\ &= 4\int_0^1 dt \int_t^1 ds\,\frac{t}{2-s} = 8\log 2 - 5. \end{aligned}$$

Summing, we get $\mathbb{E}(\Gamma^2) = I + II = \frac{5}{3}$.

It remains to prove that $\mathbb{E}(\Gamma_n) \to \mathbb{E}(\Gamma)$ and $\mathbb{E}(\Gamma_n^2) \to \mathbb{E}(\Gamma^2)$ as $n \to \infty$. By (1.15), we have

$$(5.39) \quad \mathbb{E}(C[n]) = \sum_{x:\|x\|=n} \mathbb{P}(x \in C) = [1 + o(1)]\sigma^n \times \sigma^{-(n+1)}(\sigma - 1)n\tfrac{1}{2},$$

and hence

$$(5.40) \quad \lim_{n\to\infty} \mathbb{E}(\Gamma_n) = \lim_{n\to\infty} \frac{2\sigma}{\sigma - 1}\frac{1}{n}\mathbb{E}(C[n]) = 1.$$

A similar argument applies for the second moment. Indeed, for $n_1, n_2 \to \infty$ with $n_2 \geq n_1$, we write

$$
\begin{aligned}
\mathbb{E}(C[n_1]C[n_2]) &= \sum_{x_1:\,\|x_1\|=n_1} \sum_{x_2:\,\|x_2\|=n_2} \mathbb{P}(x_1, x_2 \in C) \\
&= [1+o(1)] \sum_{k=0}^{n_1-1} \sigma^k \times \sigma^{n_1-k} \times (\sigma-1)\sigma^{n_2-k-1} \\
&\quad \times \sigma^{-[k+(n_1-k)+(n_2-k)+1]}(\sigma-1)[k+(n_1-k)+(n_2-k)] \\
&\quad \times \frac{k}{k+(n_1-k)+(n_2-k)} \frac{1}{2}\left(1 + \frac{n_1-k}{n_2} + \frac{n_2-k}{n_1}\right),
\end{aligned}
\tag{5.41}
$$

where the terms with $o$, $x_1$, $x_2$ all on a single path have been absorbed into the error term, and where we have split the sum according to the height $k$ of the most recent common ancestor of $x_1$ and $x_2$, counted the number of configurations with fixed $k$, and inserted the asymptotic formula in (1.17). For $a \in (0, 1]$, this gives

$$
\begin{aligned}
\lim_{\substack{n_1,n_2 \to \infty \\ n_1/n_2 \to a}} \mathbb{E}(\Gamma_{n_1}\Gamma_{n_2}) &= \lim_{\substack{n_1,n_2 \to \infty \\ n_1/n_2 \to a}} \left(\frac{2\sigma}{\sigma-1}\right)^2 \frac{1}{n_1 n_2} \mathbb{E}(C[n_1]C[n_2]) \\
&= 2a \int_0^1 dt\, t[1 + a(1-t) + (a^{-1} - t)] \\
&= 1 + \frac{1}{3}a(1+a),
\end{aligned}
\tag{5.42}
$$

while for $a = 0$ the limit is $2\int_0^1 dt\, t = 1$. This proves (1.26), and with $a = 1$ also $\lim_{n \to \infty} \mathbb{E}(\Gamma_n^2) = \frac{5}{3}$.

This completes the proof of (1.22).

**6. Cluster size below a given height.** In this section, we prove Theorem 1.6. The arguments and notations mirror those in Section 5. We re-use the names of the functions in Section 5 for new functions here.

6.1. *Laplace transform of $C_o[1, m]$.* For $m \in \mathbb{N}$, let $C_o[1, m]$ denote the number of vertices in the cluster of $o$ at all heights from 1 to $m$, via a fixed branch from $o$, in independent bond percolation with parameter $p \leq p_c = 1/\sigma$. For $\tau \geq 0$, let

$$f_m(p; \tau) = \mathbb{E}_p(e^{-\tau C_o[1,m]}). \tag{6.1}$$

By conditioning on the occupation status of the edge leaving the root, we see that $C_o[1, m+1]$ is 0 with probability $1-p$ and is 1 plus the sum of $\sigma$ independent copies of $C_o[1,m]$ with probability $p$. Therefore $f_m$ obeys the recursion relation

(6.2) $\qquad f_{m+1}(p;\tau) = 1 - p + pe^{-\tau}[f_m(p;\tau)]^\sigma, \qquad f_0(p;\tau) = 1.$

As before, let $g_m(p;\tau) = 1 - f_m(p;\tau)$ and $\rho = \frac{\sigma-1}{2\sigma}$. We again suppress the arguments $p$ and $\tau$. In terms of $g_m$, the recursion reads

(6.3) $\qquad\qquad\qquad g_{m+1} = F(g_m), \qquad g_0 = 0,$

with

(6.4) $\qquad\qquad F(x) = p[1 - e^{-\tau}(1-x)^\sigma], \qquad x \in [0,1].$

We will compare $g_m$ with the solution of the quadratic recursion

(6.5) $\qquad\qquad\qquad \widehat{g}_{m+1} = \widehat{F}(\widehat{g}_m), \qquad \widehat{g}_0 = 0,$

where $\widehat{F}(x)$ is the second-order approximation of $F(x)$. Thus,

(6.6) $\qquad\qquad \widehat{F}(x) = p(1 - e^{-\tau}) + \alpha x - \tfrac{1}{2}\beta x^2, \qquad x \in [0,1],$

where we abbreviate $\alpha = p\sigma e^{-\tau}$ and $\beta = (\sigma-1)\alpha$; note that $\alpha \in [0,1]$. We will compare $\widehat{g}_m$ in turn with the solution of the quadratic differential equation

(6.7) $\qquad\qquad \widetilde{g}'(t) = \widehat{F}(\widetilde{g}(t)) - \widetilde{g}(t), \qquad \widetilde{g}(0) = 0.$

The differential equation (6.7) is easily solved, as follows. By applying the linear transformation

(6.8) $\qquad\qquad\qquad y(t) = (1-\alpha) + \beta\widetilde{g}(t),$

we can write (6.7) as $y'(t) = \tfrac{1}{2}[D^2 - y(t)^2]$, where

(6.9) $\qquad\qquad\qquad D = \sqrt{(1-\alpha)^2 + 2\beta p(1 - e^{-\tau})}.$

This can be rewritten as

(6.10) $\qquad\qquad \left(\frac{1}{D+y(t)} + \frac{1}{D-y(t)}\right) y'(t) = D,$

and then integrated to give

(6.11) $\qquad \frac{D+y(t)}{D-y(t)} = \frac{D+y(0)}{D-y(0)} e^{Dt} \quad \text{or} \quad y(t) = D\frac{Ce^{Dt}-1}{Ce^{Dt}+1}$

with

(6.12) $\qquad\qquad C = \frac{D+y(0)}{D-y(0)}, \qquad y(0) = 1 - \alpha.$

Thus

(6.13) $\qquad \widetilde{g}(m) = \frac{1}{\beta}[y(m) - (1-\alpha)] = \frac{1}{\beta}\left[D\frac{Ce^{Dm}-1}{Ce^{Dm}+1} - (1-\alpha)\right].$

LEMMA 6.1. *For $m \in \mathbb{N}_0$, $p \leq p_c$ and $\tau \geq 0$,*

(6.14) $$\widetilde{g}(m) \leq g_m(p;\tau) \leq \widetilde{g}(m) + \tau + m^4 \tau^3,$$

*with $\widetilde{g}(m)$ given by (6.13), (6.9) and (6.12).*

REMARK. For $\tau \ll 1$ and $p \sim p_c$ the approximation by $\widetilde{g}(m)$ is in fact much better than stated in the lemma, though the upper bound above is sufficient for our needs.

PROOF OF LEMMA 6.1. The proof closely follows that of Lemma 5.1, but with some reversals of monotonicity. We will prove the sandwiches

(6.15) $$0 \leq g_m - \widehat{g}_m \leq m^4 \tau^3$$

and

(6.16) $$0 \leq \widehat{g}_m - \widetilde{g}(m) \leq \tau,$$

which together give the lemma.

We begin with (6.15). Since $0 \leq F'''(x) \leq \alpha \sigma^2 \leq \sigma^2$, we have

(6.17) $$0 \leq F(x) - \widehat{F}(x) \leq \tfrac{1}{6} \sigma^2 x^3, \qquad x \in [0,1].$$

Moreover, $F(0) = \widehat{F}(0) = p(1 - e^{-\tau})$, $0 \leq F'(x) \leq \alpha$, and $\widehat{F}'(x) \leq \alpha$, so

(6.18) $$0 \leq F(x) - F(0) \leq \alpha x, \qquad \widehat{F}(x) - \widehat{F}(0) \leq \alpha x.$$

As in (5.14)–(5.15), this inductively yields

(6.19) $$g_m \geq \widehat{g}_m, \qquad m \in \mathbb{N}_0.$$

In addition, $F(x) \leq x + p\tau$ and $g_0 = 0$, and therefore

(6.20) $$g_m \leq m p \tau \leq \frac{m\tau}{\sigma}, \qquad m \in \mathbb{N}_0.$$

It then follows from $F' \leq 1$, (6.17) and (6.19)–(6.20) that

(6.21) $$\begin{aligned} g_{m+1} - \widehat{g}_{m+1} &= [F(g_m) - F(\widehat{g}_m)] + [F(\widehat{g}_m) - \widehat{F}(\widehat{g}_m)] \\ &\leq (g_m - \widehat{g}_m) + \frac{1}{6} \sigma^2 g_m^3 \leq (g_m - \widehat{g}_m) + \frac{1}{6\sigma} m^3 \tau^3, \end{aligned}$$

which yields (6.15).

Next we prove (6.16). Let $\widetilde{F}(x) = h_x(1)$, where $h_x$ is the solution of

(6.22) $$h'(t) = \widehat{F}(h(t)) - h(t), \qquad h(0) = x.$$

According to (6.7), we have

(6.23) $$\widetilde{g}(m+1) = \widetilde{F}(\widetilde{g}(m)).$$

Let $x_* = [D - (1-\alpha)]/\beta \geq 0$ with $D$ defined in (6.9), and note that $\widehat{F}(x) \geq x$ for $x \in [0, x_*]$. Therefore solutions of (6.22) with $h(0) \in [0, x_*]$ are increasing on $[0, \infty)$, and the function $\widetilde{F}$ is increasing. Henceforth, we will assume the restriction $x \in [0, x_*]$. Now,

$$(6.24) \qquad h(1) - h(0) = \int_0^1 dt [\widehat{F}(h(t)) - h(t)],$$

and so, $h$ being increasing and $\widehat{F}(x) - x$ decreasing, we have

$$(6.25) \qquad \widehat{F}(h(1)) - h(1) \leq h(1) - h(0) \leq \widehat{F}(h(0)) - h(0).$$

Since $h(1) = \widetilde{F}(h(0))$, the upper bound with $h(0) = x$ gives $\widehat{F}(x) \geq \widetilde{F}(x)$. Since $\widetilde{F}$ is increasing, this inductively implies

$$(6.26) \qquad \widehat{g}_m \geq \widetilde{g}(m), \qquad m \in \mathbb{N}_0.$$

Also,

$$(6.27) \quad \widehat{g}_{k+1} - \widehat{g}_k = \widehat{F}(\widehat{g}_k) - \widehat{g}_k \leq \widehat{F}(\widetilde{g}(k)) - \widetilde{g}(k) \leq \widetilde{g}(k) - \widetilde{g}(k-1),$$

where the first inequality follows from the fact that $\widehat{F}(x) - x$ is decreasing, and the second inequality follows from the lower bound of (6.25) with $h(0) = \widetilde{g}(k-1)$ and $h(1) = \widetilde{g}(k)$.

Summing (6.27) from $k = 1$ to $k = m - 1$, we obtain $\widehat{g}_m - \widehat{g}_1 \leq \widetilde{g}(m-1) - \widetilde{g}(0)$. Since $\widetilde{g}(m-1) \leq \widetilde{g}(m)$, $\widetilde{g}(0) = g_0 = 0$ and $\widehat{g}_1 = \widehat{F}(g_0) = F(0)$, this yields (6.16). □

6.2. *Proof of Theorem* 1.6. For $k = 0, \ldots, n-1$, let $C_{k,j}[k+1, n]$ denote the contribution to the cluster size between heights $k+1$ and $n$ due to the $j$th of the $\sigma - 1$ branches emerging from the backbone at height $k$. Then $C[0, n] = n + 1 + \sum_{k=0}^{n-1} \sum_{j=1}^{\sigma-1} C_{k,j}[k+1, n]$, with the additional $n+1$ due to the backbone vertices. Conditional on $W = (W_k)_{k \in \mathbb{N}_0}$, the $C_{k,j}[k+1, n]$ are all independent. As in (5.23),

$$(6.28) \qquad \mathbb{E}(e^{-\tau C_{k,j}[k+1,n]} \mid W) = f_{n-k}(\widehat{W}_k; \tau),$$
$$k = 0, \ldots, n-1; j = 1, \ldots, \sigma - 1; \tau \geq 0,$$

with $f_m(p; \tau)$ as defined in (6.1). Consequently, since $\widehat{\Gamma}_n = \frac{1}{\rho n^2} C[0, n]$ and $f_0 = 1$,

$$\mathbb{E}(e^{-\tau \widehat{\Gamma}_n} \mid W)$$
$$(6.29) \qquad = e^{-\tau(n+1)/\rho n^2} \left( \prod_{k=0}^n f_{n-k}\left(\widehat{W}_k; \frac{\tau}{\rho n^2}\right) \right)^{\sigma-1}$$
$$= e^{-\tau(n+1)/\rho n^2} \exp\left[ (\sigma - 1) \sum_{k=0}^n \log\left(1 - g_{n-k}\left(\widehat{W}_k; \frac{\tau}{\rho n^2}\right)\right) \right].$$

Note that the prefactor $e^{-\tau(n+1)/\rho n^2}$ is equal to $1 + O(1/n)$.

Fix $\varepsilon > 0$ and $\tau \geq 0$. Since, by (6.20),

$$
(6.30) \qquad 0 \leq g_{n-k} \leq (n-k)\frac{\tau}{\rho^2 n^2} \leq O\left(\frac{\tau}{n}\right),
$$

we may linearize the logarithm and disregard the first $\varepsilon n + 1$ and the last $\varepsilon n$ terms, to get

$$
(6.31) \quad \mathbb{E}(e^{-\tau \widehat{\Gamma}_n} \mid W) = e^{-\tau(n+1)/\rho n^2} \exp\left[-\widehat{S}_n^\varepsilon(\tau, W) + O\left(\varepsilon\tau + \frac{\tau^2}{n}\right)\right]
$$

with

$$
(6.32) \qquad \widehat{S}_n^\varepsilon(\tau, W) = (\sigma - 1) \sum_{\varepsilon n < k \leq n - \varepsilon n} g_{n-k}\left(\widehat{W}_k; \frac{\tau}{\rho n^2}\right).
$$

Next, for $t \in [\varepsilon, 1 - \varepsilon]$, put $k = \lceil nt \rceil$ and $Z_n^t = n[1 - \sigma \widehat{W}_{\lceil nt \rceil}]$, and define

$$
(6.33) \quad G_n^t(x) = (\sigma - 1) n g_{\lfloor n(1-t) \rfloor}\left(\frac{1}{\sigma}\left(1 - \frac{x}{n}\right); \frac{\tau}{\rho n^2}\right) 1_{[0,n]}(x), \qquad x \in (0, \infty).
$$

Then we may write (6.32) as

$$
(6.34) \qquad S_n^\varepsilon(\tau, W) = \int_\varepsilon^{1-\varepsilon} dt\, G_n^t(Z_n^t).
$$

Now we apply Lemma 6.1, noting that $n(1 - \alpha) = x + O(\frac{\tau}{n})$, to obtain

$$
(6.35) \quad \begin{aligned} G_n^t(x) &= \left[nD\frac{Ce^{nD(1-t)} - 1}{Ce^{nD(1-t)} + 1} - x\right]\left[1 + O\left(\frac{x}{n} + \frac{\tau}{n^2}\right)\right] \\ &\quad + O\left(\frac{\tau + \tau^3}{n}\right), \qquad x \in [0, n], n \to \infty, \end{aligned}
$$

with

$$
(6.36) \qquad nD = \sqrt{x^2 + 4\tau} + O\left(\frac{\tau}{n}\right), \qquad C = \frac{nD + x}{nD - x} + O\left(\frac{\tau}{n}\right).
$$

Let $\bar{\kappa} = \bar{\kappa}(\tau, x) = \sqrt{x^2 + 4\tau}$, and put

$$
(6.37) \quad \begin{aligned} G^t(x) &= \bar{\kappa}\frac{((\bar{\kappa}+x)/(\bar{\kappa}-x))e^{(1-t)\bar{\kappa}} - 1}{((\bar{\kappa}+x)/(\bar{\kappa}-x))e^{(1-t)\bar{\kappa}} + 1} - x \\ &= \frac{4\tau}{x + \bar{\kappa}\coth[(1/2)(1-t)\bar{\kappa}]}, \qquad x \in (0, \infty). \end{aligned}
$$

Note that $G^t(x) \leq G^t(0) = 4\tau$ and that $\lim_{x \to \infty} G^t(x) = 0$ uniformly in $t \in [\varepsilon, 1-\varepsilon]$. As $n \to \infty$, $G_n^t(x)$ converges to $G^t(x)$ uniformly in $x \in (0, \infty)$ and $t \in [\varepsilon, 1-\varepsilon]$. The reasoning applied in (5.30)–(5.33) can also be applied here, to conclude that

$$S_n^\varepsilon(\tau, W) \Longrightarrow \int_\varepsilon^{1-\varepsilon} dt\, G^t(L(t)). \tag{6.38}$$

We substitute (6.38) into (6.34), and let $n \to \infty$ followed by $\varepsilon \downarrow 0$, to get

$$\lim_{n \to \infty} \mathbb{E}(e^{-\tau \widehat{\Gamma}_n}) = \mathbb{E}\left(\exp\left[-\int_0^1 dt\, G^t(L(t))\right]\right), \qquad \tau \geq 0. \tag{6.39}$$

Since $\bar{\kappa}(\tau, L(t)) = \kappa(\tau, t)$, the integral in the right-hand side equals $\widehat{S}(\tau, L)$ of (1.24), and this proves (1.23).

6.3. *First and second moment.* In this section we prove (1.25).

Taylor's expansion up to first order in $\tau$ of the integrand in (1.24) gives

$$\frac{\partial \widehat{S}}{\partial \tau}(0, L) = 2\int_0^1 dt\, \frac{1}{L(t)}[1 - e^{-(1-t)L(t)}], \tag{6.40}$$

$$-\frac{\partial^2 \widehat{S}}{\partial \tau^2}(0, L) = 8\int_0^1 dt \left(\frac{1}{2L(t)^3}[1 - e^{-2(1-t)L(t)}] - \frac{1-t}{L(t)^2}e^{-(1-t)L(t)}\right).$$

Hence, applying Lemma 3.5 for the fourth equality, we have

$$\mathbb{E}(\widehat{\Gamma}) = \mathbb{E}\left(\frac{\partial \widehat{S}}{\partial \tau}(0, L)\right) = 2\int_0^1 dt\, \mathbb{E}\left(\frac{1}{L(t)}[1 - e^{-(1-t)L(t)}]\right)$$

$$= 2\int_0^1 dt \int_t^1 ds\, \mathbb{E}(e^{-(1-s)L(t)}) = 2\int_0^1 dt \int_t^1 ds\, \frac{t}{t+1-s} = \frac{1}{2}. \tag{6.41}$$

Also, as in (5.36)–(5.38), we have $\mathbb{E}(\widehat{\Gamma}^2) = I + II$ with

$$I = \mathbb{E}\left(\left[\frac{\partial \widehat{S}}{\partial \tau}(0, L)\right]^2\right), \qquad II = \mathbb{E}\left(-\frac{\partial^2 \widehat{S}}{\partial \tau^2}(0, L)\right). \tag{6.42}$$

Using Lemma 3.5, we compute

$$I = \mathbb{E}\left(\left[2\int_0^1 dt \int_t^1 ds\, e^{-(1-s)L(t)}\right]^2\right)$$

$$= 8\int_0^1 dt_1 \int_{t_1}^1 ds_1 \int_{t_1}^1 dt_2 \int_{t_2}^1 ds_2 \frac{t_1}{t_1+(1-s_1)} \frac{t_2+(1-s_1)}{t_2+(1-s_1)+(1-s_2)}$$

$$= 8\int_0^1 dt_1 \int_{t_1}^1 ds_1 \int_{t_1}^1 dt_2 \frac{t_1}{t_1+(1-s_1)}[t_2+(1-s_1)]\log\left[\frac{2-s_1}{t_2+(1-s_1)}\right]$$

$$= \int_0^1 dt_1 \int_{t_1}^1 ds_1\, t_1[t_1+(1-s_1)]$$

(6.43)

$$\times \{4\log[t_1 + (1-s_1)] - 4\log(2-s_1) - 2\}$$

$$+ 2\int_0^1 dt_1 \int_{t_1}^1 ds_1 \frac{t_1(2-s_1)^2}{t_1+(1-s_1)}$$

$$= \int_0^1 dx \left\{ 2x^3 \log x - 2x(1+x)^2 \log(1+x) + \frac{13}{6}x^3 + \frac{8}{3}x^2 + x \right\}$$

$$= \frac{17}{8} - \frac{8}{3}\log 2,$$

where to get the next to last line we set $x = t_1 + (1-s_1)$ and interchange the two integrals. For $II$, we expand the integrand of the second line of (6.40) in powers of $L(t)$, take the expectation using that $\mathbb{E}(L(t)^k) = k!/t^k$, $k \in \mathbb{N}_0$, and sum out afterward, to obtain

(6.44)
$$II = \mathbb{E}\left(-\frac{\partial^2 \widehat{S}}{\partial \tau^2}(0, L)\right)$$
$$= \mathbb{E}\left(8\int_0^1 dt(1-t)^3 \sum_{k=0}^\infty [-(1-t)L(t)]^k \left[\frac{2^{k+2}}{(k+3)!} - \frac{1}{(k+2)!}\right]\right)$$
$$= \int_0^1 dt \left\{ 2t(2-t)^2 \log(2-t) - 2t^3 \log t - \frac{1}{2}t + \frac{1}{2}t^2 \right\}$$
$$= \frac{8}{3}\log 2 - \frac{16}{9}.$$

Summing, we get $\mathbb{E}(\widehat{\Gamma}^2) = I + II = \frac{25}{72}$.

To verify that $\mathbb{E}(\widehat{\Gamma}_n) \to \mathbb{E}(\widehat{\Gamma})$ and $\mathbb{E}(\widehat{\Gamma}_n^2) \to \mathbb{E}(\widehat{\Gamma}^2)$ as $n \to \infty$, we return to (5.39)–(5.42). Since $C[0,n] = 1 + \sum_{k=1}^n C[k]$ by definition [recall (1.2)], it follows from (5.40) that $\lim_{n\to\infty} \mathbb{E}(\widehat{\Gamma}_n) = \lim_{n\to\infty} \frac{1}{n^2} \sum_{k=1}^n k = \int_0^1 du\, u = \frac{1}{2}$. A similar calculation, based on (5.41)–(5.42), yields

(6.45)
$$\lim_{n\to\infty} \mathbb{E}(\widehat{\Gamma}_n^2) = \lim_{n\to\infty} \frac{1}{n^4} \sum_{k,l=1}^n kl \mathbb{E}(\Gamma_k \Gamma_l)$$
$$= 2\int_0^1 du \int_u^1 dv\, uv \lim_{n\to\infty} \mathbb{E}(\Gamma_{\lceil un \rceil} \Gamma_{\lceil vn \rceil})$$
$$= 2\int_0^1 du \int_u^1 dv\, uv \left[1 + \frac{1}{3}\frac{u}{v}\left(1 + \frac{u}{v}\right)\right]$$
$$= \frac{25}{18}\int_0^1 dv\, v^3 = \frac{25}{72}.$$

This completes the proof of (1.25).

**7. Proof of Theorem 1.7.** In this section we prove Theorem 1.7. The key ingredient is the following mixing property, which is proved below.

LEMMA 7.1. *Let $W^{(n)} = (W_k^{(n)})_{k \in \mathbb{N}_0}$, $n \in \mathbb{N}$, be independent realizations of $W$. There exists a sequence $(k_n)_{n \in \mathbb{N}}$ with $k_{n+1} > nk_n$, $n \in \mathbb{N}$, and a coupling of $W$ to $W^{(n)}$, $n \in \mathbb{N}$, such that with probability 1*

$$W_k = W_k^{(n)} \qquad \text{for all } k \text{ with } k_n \leq k \leq nk_n, \text{ for all but finitely many } n.$$
(7.1)

PROOF OF THEOREM 1.7. Let $C^k[n]$ denote the number of vertices in $C$ at height $n$ whose most recent ancestor on the backbone is at height at least $k$, and let $C_\infty^k[n]$ denote the same quantity for $C_\infty$. Define $\Gamma_n^k = \frac{1}{\rho n} C^k[n]$, $\Gamma_{n,\infty}^k = \frac{1}{\rho n} C_\infty^k[n]$. A small modification of the proofs of (1.4) and Theorem 1.5 shows that if $k = k(n) = o(n)$ as $n \to \infty$, then $\Gamma_n^k \Longrightarrow \Gamma$ as $n \to \infty$ under the law $\mathbb{P}$, and $\Gamma_{n,\infty}^k \Longrightarrow \Gamma_\infty$ as $n \to \infty$ under the law $\mathbb{P}_\infty$.

Since different off-backbone branches are conditionally independent given $W$, it follows from Lemma 7.1 that we may couple $C$ with independent realizations of $C^{(n)}$ such that $C^{k_n}[nk_n] = C^{(n),k_n}[nk_n]$ for all but finitely many $n$. Let

$$S_n = \frac{1}{n} \sum_{m=1}^n \Gamma_{mk_m}^{k_m}, \qquad n \in \mathbb{N}. \tag{7.2}$$

Under the law $\mathbb{P}$, all but finitely many of the summands are equal to independent random variables that converge in distribution to $\Gamma$. Consequently, $S_n \to \mathbb{E}(\Gamma)$ a.s. as $n \to \infty$ under $\mathbb{P}$. On the other hand, let

$$S_{n,\infty} = \frac{1}{n} \sum_{m=1}^n \Gamma_{mk_m,\infty}^{k_m}, \qquad n \in \mathbb{N}. \tag{7.3}$$

Then, under the law $\mathbb{P}_\infty$, the summands are already independent and converge in distribution to $\Gamma_\infty$, so that $S_{n,\infty} \to \mathbb{E}_\infty(\Gamma_\infty)$ a.s. as $n \to \infty$ under $\mathbb{P}$. Since $\mathbb{E}(\Gamma) \neq \mathbb{E}_\infty(\Gamma_\infty)$, it follows that the random variables $S_n$ and $S_{n,\infty}$ (which are actually the *same* variables under different laws) a.s. have unequal limits as $n \to \infty$. Therefore, the laws of IPC and IIC are mutually singular, since the IPC is supported on clusters for which the limit is $\mathbb{E}(\Gamma)$, whereas the IIC is supported on clusters which have the different limit $\mathbb{E}_\infty(\Gamma_\infty)$. □

PROOF OF LEMMA 7.1. Given a realization of the Markov chain $W$, consider a realization $W'$ that uses the same sequence of random variables $X_k$ [recall (3.13)]. Since $R(\cdot)$ is increasing, it is possible to arrange the coupling so that whenever the process with the lower value jumps, the process with the higher value also jumps (necessarily to the same $X_k$). Let $\tau$ be the first time such a jump occurs, and note that $W_k = W'_k$ for all $k > \tau$.

Suppose $W'_k > W_k$ for some $k$. This inequality is preserved until some step at which $W'$ jumps to a value in $(p_c, W_k)$. This happens with probability $\theta(W_k)R(W'_k)$, while $W$ jumps with probability $\theta(W_k)R(W_k)$. Thus the probability that the processes coalesce at such a jump is $R(W_k)/R(W'_k) \geq \rho = \frac{\sigma-1}{2\sigma} \geq 1/4$ [recall (3.10)].

Let $Z_k = W_k \wedge W'_k$. The above shows that whenever $Z$ decreases there is probability at least $\rho$ that $W$ and $W'$ coalesce. Since $\lim_{k \to \infty} Z_k = p_c$, $Z$ decreases infinitely often and the processes coalesce at some finite time.

We can now construct the desired coupling. Let $W^{(n)}$ be independent realizations of $W$, derived with random variables $X_k^{(n)}$. We fix $k_1 = 1$ arbitrarily, and define $k_2, k_3, \ldots$ inductively as follows. Having chosen $k_{n-1}$, consider the above coupling of $W$ with $W' = W^{(n)}$, started at time $(n-1)k_{n-1}$. Note that at $k = (n-1)k_{n-1}$, the values $W_k$ and $W_k^{(n)}$ are independent. By the above argument, there is some a.s. finite $\tau$ so that $W_k = W_k^{(n)}$ for any $k > \tau$. We can, therefore, select $k_n > (n-1)k_{n-1}$ so that $\mathbb{P}(\tau \geq k_n) < n^{-2}$, and hence $\mathbb{P}(W_{k_n} \neq W_{k_n}^{(n)}) < n^{-2}$. By the Borel–Cantelli lemma, we have $W_{k_n} = W_{k_n}^{(n)}$ for all but finite many $n$ a.s., and this equality holds up to $nk_n$. □

**8. Proof of Theorems 1.8 and 1.9.** In this section, we prove Theorem 1.9. As mentioned in Section 1.3.5, Theorem 1.8 then follows via [1, Example 1.9(ii)]. We retain the terminology of Lemma 6.1 and its proof, and begin with two technical estimates.

LEMMA 8.1. *There is a $c > 0$ such that, for $p = \frac{1}{\sigma}(1-\delta)$,*

$$D \geq [1 + o(1)]\sqrt{\tau}, \qquad \delta, \tau \downarrow 0, \tag{8.1}$$

$$D - (1-\alpha) \geq c\left(\frac{\tau}{\delta} \wedge \sqrt{\tau}\right), \qquad 0 < \delta, \tau \ll 1. \tag{8.2}$$

PROOF. Recall the definition of $\alpha, \beta$ below (6.6). As $\delta, \tau \downarrow 0$,

$$\begin{aligned}
1 - \alpha &= 1 - (1-\delta)e^{-\tau} \sim \delta + \tau, \\
2\beta p(1 - e^{-\tau}) &= 4\rho(1-\delta)^2 e^{-\tau}(1 - e^{-\tau}) \sim 4\rho\tau.
\end{aligned} \tag{8.3}$$

For (8.1), we use $4\rho \geq 1$ [recall (1.3) with $\sigma \geq 2$] to obtain

$$D = \sqrt{(1-\alpha)^2 + 2\beta p(1-e^{-\tau})} \geq [1+o(1)]\sqrt{\tau}, \qquad \delta, \tau \downarrow 0. \tag{8.4}$$

For (8.2), note that

$$D \leq 2[(\delta + \tau) \vee \sqrt{4\rho\tau}], \qquad 0 < \delta, \tau \ll 1, \tag{8.5}$$

to obtain

$$D - (1 - \alpha) = \frac{D^2 - (1-\alpha)^2}{D + (1-\alpha)}$$

(8.6)
$$\geq \frac{2\beta p(1 - e^{-\tau})}{2D} \geq \frac{\rho\tau}{D} \geq c\left(\frac{\tau}{\delta} \wedge \sqrt{\tau}\right), \qquad 0 < \delta, \tau \ll 1.$$

□

LEMMA 8.2. *There is a $c > 0$ such that, for $n \gg 1$ and $\tau \geq n^{-2}$,*

(8.7) $$g_n\left(\frac{1}{\sigma}(1-\delta); \tau\right) \geq c[D - (1-\alpha)].$$

PROOF. We use Lemma 6.1 and the inequalities $\beta < \sigma$ and $C > 1$ [recall (6.12)], to estimate

(8.8)
$$\begin{aligned}
g_n\left(\frac{1}{\sigma}(1-\delta); \tau\right) &\geq \frac{1}{\beta}\left[D\frac{Ce^{Dn} - 1}{Ce^{Dn} + 1} - (1-\alpha)\right] \\
&= \frac{1}{\beta}\left[D - (1-\alpha) - \frac{2D}{Ce^{Dn} + 1}\right] \\
&\geq \frac{D - (1-\alpha)}{\sigma}\left[1 - \frac{2D}{[D + (1-\alpha)]e^{Dn}}\right] \\
&\geq \frac{D - (1-\alpha)}{\sigma}[1 - 2e^{-Dn}] \geq c[D - (1-\alpha)],
\end{aligned}$$

where in the last inequality we use (8.1), which implies that $Dn \geq [1 + o(1)]n\sqrt{\tau} \geq [1 + o(1)]$. □

PROOF OF THEOREM 1.9. For convenience we restrict ourselves to $n$ divisible by 3, which suffices. Our goal is to get a bound uniform in $n$ for

(8.9) $$9n^2 \mathbb{E}\left(\frac{1}{C[0, 3n]}\right) = 9n^2 \int_0^\infty d\tau \, \mathbb{E}(e^{-\tau C[0, 3n]}),$$

where the equality follows from Fubini's theorem.

The contribution to the right-hand side of (8.9) due to $0 \leq \tau < n^{-2}$ is bounded by 9. Moreover, due to the presence of the backbone, we have $C[0, 3n] \geq 3n + 1$, and so

(8.10) $$\int_{n^{-1}\log n}^\infty d\tau \, \mathbb{E}(e^{-\tau C[0,3n]}) \leq \int_{n^{-1}\log n}^\infty d\tau \, e^{-3n\tau} = \frac{1}{3n^4}.$$

Thus, it suffices to show that

(8.11) $$\sup_{n \in \mathbb{N}} n^2 \int_{n^{-2}}^{n^{-1}\log n} d\tau \, \mathbb{E}(e^{-\tau C[0,3n]}) < \infty.$$

As in (6.29), we have

$$\mathbb{E}(e^{-\tau C[0,3n]} \mid W) \leq \prod_{k=0}^{3n}(1 - g_{3n-k}(\widehat{W}_k;\tau))^{\sigma-1}$$

(8.12)
$$\leq \exp\left[-\sum_{k=n}^{2n-1} g_{3n-k}(\widehat{W}_k;\tau)\right],$$

where the second inequality uses $g_m \geq 0$ and $\sigma \geq 2$. We know that $g_m(p;\tau)$ is increasing in $p$ and $m$ [because $f_m(p;\tau)$ in (6.1) is decreasing in $p, m$]. Since $\widehat{W}_k$ is increasing in $k$, we may therefore estimate

(8.13) $$\sum_{k=n}^{2n-1} g_{3n-k}(\widehat{W}_k;\tau) \geq n g_n(\widehat{W}_n;\tau).$$

Thus, to get (8.11), it suffices to show that

(8.14) $$\sup_{n \in \mathbb{N}} n^2 \int_{n^{-2}}^{n^{-1}\log n} d\tau \mathbb{E}(\exp[-n g_n(\widehat{W}_n;\tau)]) < \infty.$$

Let $\delta_n = 1 - \sigma \widehat{W}_n$. By Lemma 3.2, $\mathbb{P}(\delta_n \geq 1/\sqrt{n}) \leq \exp[-c'\sqrt{n}]$ for some $c' > 0$, and so we may restrict the integral in (8.14) to the event $\{\delta_n < 1/\sqrt{n}\}$. By Lemmas 8.1–8.2, there is a $c > 0$ such that, for $n$ sufficiently large to give $\delta_n, \tau \ll 1$,

(8.15)
$$\exp[-n g_n(\widehat{W}_n;\tau)] = \exp\left[-n g_n\left(\frac{1}{\sigma}(1-\delta_n);\tau\right)\right]$$
$$\leq \exp\left[-cn\left(\frac{\tau}{\delta_n} \wedge \sqrt{\tau}\right)\right] \leq e^{-cn\tau/\delta_n} + e^{-cn\sqrt{\tau}}.$$

Thus, to get (8.14), it suffices to show that

(8.16)
$$\sup_{n \in \mathbb{N}} n^2 \int_0^\infty d\tau [\mathbb{E}(e^{-cn\tau/\delta_n}) + e^{-cn\sqrt{\tau}}]$$
$$= \frac{1}{c} \sup_{n \in \mathbb{N}} \mathbb{E}(n\delta_n) + \int_0^\infty dt e^{-c\sqrt{t}} < \infty.$$

But, by Lemma 3.2, the last supremum is finite, and so the proof is complete. □

REMARK. Note that the IIC corresponds to $\delta_n = 0$ for all $n$, and that it follows from [2], equation (2.17) that the term $\int_0^\infty dt e^{-c\sqrt{t}}$ in (8.16) is an upper bound for the corresponding IIC expectation $\mathbb{E}_\infty(n^2/C_\infty[0,n])$. The additional term in (8.16) is a reflection of the fact that the IPC is smaller than the IIC, and the uniform boundedness of the expectation $\mathbb{E}(n\delta_n)$ is consistent with the scaling of the $\widehat{W}$-process in Corollary 3.4.

**9. Proofs for the incipient infinite cluster.** In this section, we give quick proofs of the statements made in Section 1.2. First we look at the structural representation of $C_\infty$ under the law $\mathbb{P}_\infty$, and then we turn to the $r$-point functions and to the cluster size at and below a given height.

*Structural representation.* Let $\{o \to n\}$ denote the event that the root is connected to a vertex that is distance $n$ from the root. The law $\mathbb{P}_\infty$ is defined as the limit

$$(9.1) \qquad \mathbb{P}_\infty(E) = \lim_{n \to \infty} \mathbb{P}_{p_c}(E \mid o \to n) \qquad \forall \text{ cylinder event } E.$$

Let $a_1, \ldots, a_\sigma$ denote the neighbors of $o$. Then

$$(9.2) \qquad \mathbb{P}_{p_c}(E \mid o \to n) = \frac{1}{\sigma} \sum_{1 \leq i \leq \sigma} \mathbb{P}_{p_c}(E \mid o \to a_i, a_i \xrightarrow{o} n) + O\left(\frac{1}{n}\right),$$

where $\xrightarrow{o}$ means a connection avoiding $o$, and the error term covers the case of two or more disjoint connections from $o$ to $n$. Suppose that $E$ is an elementary event, that is, $E$ determines the state of all the edges it depends on. Let $\mathcal{T}_i$ denote the $i$th branch of $\mathcal{T}_\sigma$ from $o$, $E_i$ the restriction of $E$ to $\mathcal{T}_i$, and $\mathbb{P}^i_{p_c}$ the critical percolation measure on $\mathcal{T}_i$. Then $E = \otimes_{i=1}^\sigma E_i$ and

$$(9.3) \quad \mathbb{P}_{p_c}(E \mid o \to a_i, a_i \xrightarrow{o} n) = \left[\prod_{\substack{1 \leq j \leq \sigma \\ j \neq i}} \mathbb{P}^j_{p_c}(E_j)\right] \mathbb{P}^i_{p_c}(E_i \mid 0 \to a_i \to n).$$

Now let $n \to \infty$ and use that the last factor tends to $\mathbb{P}^i_\infty(E_i)$, which is defined as the law on $\mathcal{T}_i$ under which the edge between $o$ and $a_i$ is open and from $a_i$ there is an IIC. Then, using (9.2), we obtain

$$(9.4) \qquad \mathbb{P}_\infty(E) = \frac{1}{\sigma} \sum_{1 \leq i \leq \sigma} \left[\prod_{\substack{1 \leq j \leq \sigma \\ j \neq i}} \mathbb{P}^j_{p_c}(E_j)\right] \mathbb{P}^i_\infty(E_i).$$

What this equation says is that the IIC can be grown recursively by opening a random edge, putting critical percolation clusters in the branches emerging from the bottom of this edge, and proceeding to grow from the top of this edge. □

PROOF OF (1.1). Pick $y \in \partial \mathcal{S}(\vec{x})$ and recall that $B_y$ is the event that $y$ is in the backbone. We have

$$(9.5) \qquad \mathbb{P}_\infty(B_y) = \sigma^{-\|y\|},$$

$$\mathbb{P}_\infty(\vec{x} \in C_\infty \mid B_y) = \sigma^{-(N-\|y\|+1)}.$$

Indeed, the first line comes from the fact that the backbone is uniformly random, while the second line uses that $\mathcal{S}(\vec{x})$ has $N - \|y\| + 1$ edges off the path from $o$ to $y$. We multiply the two equations in (9.5), sum over $y$, and use that the cardinality of $\partial \mathcal{S}(\vec{x})$ is $N(\sigma-1)+\sigma$, to get the first line of (1.1). Then we divide the product by the sum, to get the second line of (1.1). □

PROOF OF (1.4)–(1.5). In view of Proposition 2.1, Lemma 2.2 and Proposition 3.3, the IIC corresponds to taking the limit when $L \equiv 0$. In this case $S(\tau, L)$ in (1.21) reduces to

$$(9.6) \qquad S_\infty(\tau) = 2\tau \int_0^1 dt \frac{1}{1+\tau(1-t)} = 2\log(1+\tau).$$

Hence we get

$$(9.7) \qquad \mathbb{E}_\infty(e^{-\tau\Gamma}) = e^{-S_\infty(\tau)} = (1+\tau)^{-2}.$$

Similarly, $\widehat{S}(\tau, L)$ in (1.24) reduces to

$$(9.8) \qquad \widehat{S}(\tau, L) = 2\sqrt{\tau} \int_0^1 dt \tanh((1-t)\sqrt{\tau}) = 2\log\cosh(\sqrt{\tau}).$$

Hence we get

$$(9.9) \qquad \mathbb{E}_\infty(e^{-\tau\widehat{\Gamma}}) = [\cosh(\sqrt{\tau})]^{-2}.$$

A more formal proof of (1.4)–(1.5) can be obtained along the lines of Sections 5–6. This requires that we set $p = p_c = \frac{1}{\sigma}$ in Lemmas 5.1 and 6.1 and repeat the estimates in Sections 5.2 and 6.2, which in fact simplify considerably. □

**Acknowledgments.** OA held a PIMS postdoctoral fellowship. JG holds an NSERC Canada Graduate Scholarship. FdH is grateful to PIMS and the Department of Mathematics of UBC for hospitality during a sabbatical leave from January to August 2006.

O. Angel  
Department of Statistics  
University of Toronto  
Toronto, Ontario  
Canada M5S 2E4  
E-mail: angel@utstat.toronto.edu

J. Goodman  
G. Slade  
Department of Mathematics  
University of British Columbia  
Vancouver, British Columbia  
Canada V6T 1Z2  
E-mail: jgoodman@math.ubc.ca  
slade@math.ubc.ca

F. den Hollander  
Mathematical Institute  
Leiden University  
P.O. Box 9512  
2300 RA Leiden  
The Netherlands  
and  
EURANDOM  
P.O. Box 513  
5600 MB Eindhoven  
The Netherlands  
E-mail: denholla@math.leidenuniv.nl